\documentclass{article}
\usepackage{amssymb}
\usepackage{amsmath,color}
\usepackage[hidelinks]{hyperref}

\usepackage[square,numbers]{natbib}         
\newtheorem{lem}{Lemma}[section]
\newtheorem{cor}[lem]{Corollary}
\newtheorem{teo}[lem]{Theorem}
\newtheorem{os}[lem]{Remark}

\newtheorem{prop}[lem]{Proposition}
\newtheorem{esem}[lem]{Example}

\newcommand{\qed}{\thinspace\null\nobreak\hfill\hbox{\vbox{\kern-.2pt\hrule
 height.2pt depth.2pt\kern-.2pt\kern-.2pt \hbox to2.5mm{\kern-.2pt\vrule
 width.4pt \kern-.2pt\raise2.5mm\vbox to.2pt{}\lower0pt\vtop
 to.2pt{}\hfil\kern-.2pt \vrule
 width.4pt \kern-.2pt}\kern-.2pt\kern-.2pt\hrule height.2pt depth.2pt
 \kern-.2pt}}\par\medbreak}

\newcommand{\R}{\mathbb{R}}

\newcommand{\C}{\mathbb{C}}

\newcommand{\N}{\mathbb{N}}
\newcommand{\Z}{\mathbb{Z}}

\newcommand{\Rp}{\textrm{\emph{Re}\,}}

\newcommand{\eps}{\varepsilon}

\newcommand{\ds}{\displaystyle}

\date{}
\begin{document}

\title{ A unified approach to degenerate problems in the half-space}
\author{G. Metafune \thanks{Dipartimento di Matematica e Fisica ``Ennio De Giorgi'', Universit\`a del Salento, C.P.193, 73100, Lecce, Italy.
-mail:  giorgio.metafune@unisalento.it}\qquad L. Negro \thanks{Dipartimento di Matematica e Fisica  ``Ennio De
Giorgi'', Universit\`a del Salento, C.P.193, 73100, Lecce, Italy. email: luigi.negro@unisalento.it} \qquad C. Spina \thanks{Dipartimento di Matematica e Fisica``Ennio De Giorgi'', Universit\`a del Salento, C.P.193, 73100, Lecce, Italy.
e-mail:  chiara.spina@unisalento.it}}

\maketitle
\begin{abstract}
\noindent 
We study elliptic and parabolic problems governed by the  singular elliptic   operators 
 \begin{equation*}
\mathcal L =y^{\alpha_1}\Delta_{x} +y^{\alpha_2}\left(D_{yy}+\frac{c}{y}D_y  -\frac{b}{y^2}\right),\qquad\alpha_1, \alpha_2 \in\R
\end{equation*}
in the half-space $\R^{N+1}_+=\{(x,y): x \in \R^N, y>0\}$.
 
\par

\end{abstract}

\section{Introduction}
 In this paper we study solvability and regularity of elliptic and parabolic problems associated to the  degenerate   operators 
\begin{equation*} \label{defL}
\mathcal L =y^{\alpha_1}\Delta_{x} +y^{\alpha_2}\left(D_{yy}+\frac{c}{y}D_y  -\frac{b}{y^2}\right) \quad {\rm and}\quad D_t- \mathcal L
\end{equation*}
in the half-space $\R^{N+1}_+=\{(x,y): x \in \R^N, y>0\}$ or  in $(0, \infty) \times \R^{N+1}_+$.

Here $b,\ c$ are constant real coefficients and we use 
$
L_y=D_{yy}+\frac{c}{y}D_y  -\frac{b}{y^2}.
$
 Note that singularities in the lower order terms appear when either $b$ or $c$ is different from 0. When $b=0$, then $L_y$ is a Bessel operator and we shall denote it by $B_y$.

The real numbers $\alpha_1, \alpha_2$ satisfy $\alpha_2<2$ and $\alpha_2-\alpha_1<2$ but are not assumed to be nonnegative. The reasons for these restrictions will be explained later in this introduction.

$\mathcal L$ is the sum of a degenerate diffusion $y^{\alpha_1} \Delta_x$, tangential to $\partial \R^{N+1}_+$, and of a 1d degenerate normal diffusion $y^{\alpha_2}L_y$ which commute only when $\alpha_1=0$.
 It satisfies the scaling property $$I_s^{-1} \mathcal L I_s=s^{2-{\alpha_2}} \mathcal L, \quad I_su(x,y)=u(s^{1-\frac{(\alpha_2-\alpha_1)}{2}}x, sy).$$

When $\alpha_1=\alpha_2=0$, $\mathcal L=\Delta_{x} +\left(D_{yy}+\frac{c}{y}D_y  -\frac{b}{y^2}\right)$ reduces to the so-called Caffarelli-Silvestre extension operators, studied in detail in \cite{MNS-Caffarelli} under Dirichlet and Neumann boundary conditions. We refer the reader also to  \cite{dong2020parabolic}, \cite{dong2021weighted} for the case $b=0$  and with variable coefficients.

The case $\alpha_1=0$ and $\alpha_2$ arbitrary (even without the restriction $\alpha_2<2$) can be easily reduced to that above, using the methods of this paper and, in particular, the transformation of Section \ref{Section Degenerate}.

The case $\alpha_1=\alpha_2=1$ and $b=0$,  namely $ \mathcal L =y\left (\Delta_{x} +D_{yy}\right)+cD_y $,  is also widely treated in the  literature on degenerate problems. A comparison of our results with those already known is done in section 7.4

When $\alpha_2=0$ our operators generalize the class of Baouendi-Grushin operators $\mathcal L =y^{\alpha}\Delta_{x} +D_{yy}$, to which they reduce when $c=b=0$. A comparison with known results is done in Section 7.3 and here we only point out that we allow also negative $\alpha$.

Finally, we mention that kernel estimates for operators in divergence form and with normal and tangential degeneracy on the hyperplane $\{y=0\}$ have been obtained in \cite{Robinson-Sikora2008, Robinson-Sikora2014}.

\medskip
The aim of this paper is to provide a unified approach which allows to prove elliptic and parabolic $L^p$ estimates and solvability of the associated problems.
In the language of semigroup theory, we prove that $\mathcal L$ generates an analytic semigroup, characterize its domain as a weighted Sobolev space and show that it has maximal regularity, which means that both $D_t v$ and $\mathcal L v$ have the same regularity as $(D_t -\mathcal L) v$.

Surprisingly enough, the case $\alpha_1=\alpha_2$ implies all other cases by a change of variables, as described in Section \ref{Section Degenerate}. However this modifies the underlying measure and the procedure works if one is able to deal with the simpler case $\mathcal L=y^{\alpha}(\Delta_x+L_y)$ in all the scale of $L^p_m$ spaces, where $L^p_m=L^p(\R^{N+1}_+; y^m dxdy)$. A similar simplification holds also for the 1d operator $L_y$: it is sufficient to deal in full generality with the case where $b=0$, that is when $L_y$ is  a Bessel operator. Finally, it is sufficient to deal only with Neumann boundary conditions, since the case of Dirichlet boundary conditions is again deduced by a change of variables.


The operators $\mathcal L$, $ D_t-\mathcal L$  are studied through estimates like 
\begin{equation} \label{closedness}
\|y^\alpha \Delta_x u\|_{p,m}+\|y^\alpha B_y u\|_{p,m} \le C\| \mathcal L u\|_{p,m}, \quad 
\end{equation}
and
\begin{equation} \label{closedness1}
\|D_t u\|_{p,m}+\|\mathcal Lu\|_{p,m} \le C\| (D_t-\mathcal L) u\|_{p,m}, \quad 
\end{equation}
where the $L^p$ norms are taken over $\R_+^{N+1}$ and on $(0, \infty) \times \R_+^{N+1}$ respectively. This kind of estimates are quite natural in this context but not easy to prove. Of course they imply $\|y^\alpha D_{x_ix_j}u\|_{p,m} \le C\| \mathcal Lu\|_{p,m}$, by the Calder\'{o}n-Zygmund inequalities in the $x$-variables, and can be restated by saying that $\mathcal L$ is closed on $D(y^\alpha\Delta_x) \cap D(y^\alpha B_y)$ or that $y^\alpha \Delta_x  \mathcal L^{-1}$ is bounded. Note that the weaker inequality \eqref{closedness} with $\| \mathcal Lu\|_{p,m}+\|u\|_{p,m}$ on the right hand side implies the homogeneous one, as stated, by scaling.

Let us explain how to obtain \eqref{closedness}. Assuming that $y^\alpha(\Delta_x u+B_yu)=f$ and taking the Fourier transform with respect to $x$ (with covariable $\xi$) we obtain $-|\xi|^2 \hat u(\xi,y)+B_y \hat u(\xi,y)=y^{-\alpha}\hat f(\xi,y)$ and then $y^\alpha|\xi|^2 \hat u(\xi,y)=-y^\alpha |\xi|^2 (|\xi|^2-B_y)^{-1}y^{-\alpha}\hat f (\xi,y)$. 
Denoting by ${\cal F}$ the Fourier transform with respect to $x$  we get 
$$
y^\alpha\Delta_x  \mathcal  L^{-1}=-{\cal F}^{-1}\left (y^\alpha|\xi|^2(|\xi|^2-B_y)^{-1} y^{-\alpha} \right) {\cal F}
$$
and the boundedness of $y^\alpha\Delta_x  \mathcal L^{-1}$ is equivalent to that of the multiplier $\xi \in \R^N \to M(\xi)=y^\alpha|\xi|^2(|\xi|^2-B_y)^{-1}y^{-\alpha} $ in $L^p(\R^N; L^p_m(0,\infty))=L^p_m(\R_+^{N+1})$. 

We prove this by a vector valued Mikhlin multiplier theorem which rests on square function estimates for the family $M(\xi)$ and its derivatives. The strategy for proving  \eqref{closedness1} is similar after taking the Fourier transform with respect to $t$.

Both the elliptic and parabolic estimates above  share the name ``maximal regularity" even though this term is often restricted to the parabolic case. 
We refer to  \cite{KunstWeis} and the new books \cite{WeisBook1}, \cite{WeisBook2} for 
the functional analytic approach to  maximal regularity we use. The whole theory relies on a deep interplay between harmonic analysis and structure theory of Banach spaces but largely simplifies when the underlying Banach spaces are $L^p$ spaces, by the use of classical square function estimates. This last approach has been employed extensively in \cite{DenkHieberPruss}, showing that uniformly parabolic operators have maximal regularity, under very general boundary conditions.

\medskip

However the a-priori estimates  \eqref{closedness} and \eqref{closedness1} are not sufficient for the solvability of the equation $\lambda u -\mathcal L u=f$. In fact, $\mathcal L$ is not dissipative unless additional restrictions on the parameters and on the underlying measure are assumed, see Section 9.1, and approximation methods with uniformly parabolic operators do not need to converge.

In order to prove existence results, or generation results in the language of semigroups,  we use that 
the operator valued map $\xi \in \R^N \to N(\xi)=(\lambda+y^\alpha|\xi|^2-y^\alpha B_y)^{-1}$, $\lambda\in \C_+$, is a Fourier multiplier in $L^p(\R^N; L^p_m(0,\infty))=L^p_m(\R_+^{N+1})$, see \cite[Section 8]{MNS-PerturbedBessel} where the relevant  one dimensional degenerate operators are studied in detail. 

\medskip

Before describing the content of the sections, let us explain the meaning of the restrictions $\alpha_2<2$, $\alpha_2-\alpha_1<2$. 

Let us first consider the case where $\alpha_1=\alpha_2=\alpha$, so that the unique requirement is $\alpha<2$. It turns out that when $\alpha \geq 2$ the problem is easily treated in the strip $\R^N\times [0,1]$ in the case of the Lebesgue measure, see \cite{FornMetPallScn5}, and all problems are due to the strong diffusion at infinity. The case $\alpha \geq 2$ in the strip $\R^N \times [1, \infty[$ requires therefore new investigation even though the operator $y^\alpha L_y$ alone can be treated for any $\alpha \in \R$, by the similarity transformation of Section 3.

When $\alpha_1 \neq \alpha_2$, the change of variables of Section 3, namely $T_{0, \frac{\alpha_1-\alpha_2}{2}}$, transforms $y^{\alpha_1}\Delta_x+y^{\alpha_2}B_y$ into $y^{\alpha}(\Delta_x+\tilde {B}_y)$, $\alpha=\frac{2\alpha_1}{\alpha_1-\alpha_2+2}$. However, the strip $\R^N\times [0,1]$ is mapped into itself only when $\alpha_2-\alpha_1<2$.
Under this condition it is possible, though not treated in this paper, that the restriction $\alpha_2<2$ can be removed, at least when the operator is studied in $\R^N\times [0,1]$ rather than in $\R^{N+1}_+$. But dealing with the case $\alpha_2-\alpha_1 \geq 2$ requires further investigation, as explained above. 

Assuming, in addition, that $\alpha_1 \geq 0$, the range of parameters for which we prove solvability is optimal, since it coincides with that of $L_y$. However, when $\alpha_1<0$ it can happen that $\mathcal L$ generates in a range of parameters for which the domain is less regular. We discuss these phenomena in Section 9.2, see in particular Example \ref{noreg}, without demanding for completeness. 

\medskip

Most of the result of this paper can be extended to  operators with variable coefficients
$$
\mathcal L=y^{\alpha_1}\sum_{i,j=1}^N a_{ij}(t,x,y)D_{x_ix_j} +y^{\alpha_2}\left(D_{yy}+\frac{c}{y}D_y  -\frac{b}{y^2}\right) 
$$
assuming uniform ellipticity and appropriate continuity of the matrix $(a_{ij})$. In fact, the case of constant coefficients $(a_{ij})$ follows by a linear change of the $x$-variables and allows to use perturbation methods. The situation is easier in a finite strip $\R^N \times [0,1]$ and for positive $\alpha_1, \alpha_2$ since the powers $y^{\alpha_1}, y^{\alpha_2}$ are bounded. First order terms like $y^{\frac{\alpha_1}{2}}b(t,x,y) \cdot \nabla_x$ with $b$ bounded can be also added by perturbation. We shall deal with these consequences in a subsequent paper in order to treat operators which degenerate near the boundary of a domain with (possibly) different rates along normal and tangential directions.

\medskip

The paper is organized as  follows. In Section 2 we briefly recall the harmonic analysis background needed in the paper,  as  square function estimates, $\mathcal R$-boundedness  and a vector valued multiplier theorem. 

In Section 3, we exploit an elementary change of variables, in a functional analytic setting, to reduce our operators to the simpler case where $\alpha_1=\alpha_2$.



In Sections 4 and  5 we recall some preliminary results concerning anisotropic weighted Sobolev spaces and  one-dimensional Bessel operators.

In Section 6, which is the core of the paper,  we prove generation results, maximal regularity and domain characterization for the operator $y^\alpha\Delta_x+y^\alpha B_y^n$ where $B_y^n$ is the Bessel operator with Neumann boundary conditions. The general case both with Neumann and Dirichlet boundary conditions will be deduced in Sections 7 and 8  using the isometry of Section 3.

In Section 9, we complement our results by characterizing the contractivity range and investigating uniqueness. We show that many examples of degenerate operators, like the  Baouendi-Grushin operators, are special cases of ours and that our results improve those already existing in the literature.

\bigskip
\noindent\textbf{Notation.} For $N \ge 0$, $\R^{N+1}_+=\{(x,y): x \in \R^N, y>0\}$. For $m \in \R$ we consider the measure $y^m dx dy $ in $\R^{N+1}_+$ and  we write $L^p_m(\R_+^{N+1})$ for  $L^p(\R_+^{N+1}; y^m dx dy)$ and often only $L^p_m$ when $\R^{N+1}_+$ is understood. Similarly $W^{k,p}_m(\R^{N+1}_+)=\{u \in L^p_m(\R^{N+1}_+): \partial^\alpha u \in  L^p_m(\R^{N+1}_+) \quad |\alpha| \le k\}$. 

$\C^+=\{ \lambda \in \C: \Rp \lambda >0 \}$ and, for $|\theta| \leq \pi$, we denote by  $\Sigma_{\theta}$  the open sector $\{\lambda \in \C: \lambda \neq 0, \ |Arg (\lambda)| <\theta\}$.
\bigskip

\section{Vector-valued harmonic analysis}
Regularity properties for $\mathcal L=y^{\alpha_1}\Delta_x+y^{\alpha_2}L_y$ follow once we prove the estimate
\begin{equation}  \label{regularity}
\|y^{\alpha_1}\Delta_x u\|_p+\|y^{\alpha_2}L_y u\|_p \le C \| \mathcal L u\|_p
\end{equation}
 where the $L^p$ norms are taken over $\R^{N+1}_+$ on a sufficiently large set of functions $u$. This is equivalent to saying that the domain of $\mathcal L $ is the intersection of the domain of $y^{\alpha_1}\Delta_x$ and $y^{\alpha_2}L_y$ (after appropriate tensorization) or that the operator $y^{\alpha_1}\Delta_x {\mathcal L}^{-1}$ is bounded. This strategy arose first in the study of maximal regularity of parabolic problems, that is for the equation $u_t=Au+f, u(0)=0$ where $A$ is the generator of an analytic semigroup on a Banach space $X$. Estimates like
$$
\|u_t\|_p+\|Au\|_p \le \|f\|_p
$$
where now the $L^p$ norm is that of $L^p([0, T[;X)$ can be interpreted as closedness of $D_t-A$ on the intersection of the respective domains or, equivalently, boundedness of the operator $A(D_t-A)^{-1}$ in $L^p([0, T[;X)$.

Nowadays this strategy is well established and relies on Mikhlin vector-valued multiplier theorems.
Let us state the relevant definitions and main results we need, referring the reader to \cite{DenkHieberPruss}, \cite{Pruss-Simonett} or \cite{KunstWeis}.

Let ${\cal S}$ be a subset of $B(X)$, the space of all bounded linear operators on a Banach space $X$. ${\cal S}$ is $\mathcal R$-bounded if there is a constant $C$ such that
$$
\|\sum_i \eps_i S_i x_i\|_{L^p(\Omega; X)} \le C\|\sum_i \eps_i  x_i\|_{L^p(\Omega; X)} 
$$
for every finite sum as above, where $(x_i ) \subset X, (S_i) \subset {\cal S}$ and $\eps_i:\Omega \to \{-1,1\}$ are independent and symmetric random variables on a probability space $\Omega$. The smallest constant $C$ for which the above definition holds is the $\mathcal R$-bound of $\mathcal S$,  denoted by $\mathcal R(\mathcal S)$.
It is well-known that this definition does not  depend on $1 \le p<\infty$ (however, the constant $\mathcal R(\mathcal S)$ does) and that $\mathcal R$-boundedness is equivalent to boundedness when $X$ is an Hilbert space.
When $X$ is an $L^p$ space (with respect to any $\sigma$-finite measure), testing  $\mathcal R$-boundedness is equivalent to proving square functions estimates, see \cite[Remark 2.9 ]{KunstWeis}.

\begin{prop}\label{Square funct R-bound} Let ${\cal S} \subset B(L^p(\Sigma))$, $1<p<\infty$. Then ${\cal S}$ is $\mathcal R$-bounded if and only if there is a constant $C>0$ such that for every finite family $(f_i)\in L^p(\Sigma), (S_i) \in {\cal S}$
$$
\left\|\left (\sum_i |S_if_i|^2\right )^{\frac{1}{2}}\right\|_{L^p(\Sigma)} \le C\left\|\left (\sum_i |f_i|^2\right)^{\frac{1}{2}}\right\|_{L^p(\Sigma)}.
$$
\end{prop}
The best constant $C$ for which the above square functions estimates hold satisfies $\kappa^{-1} C \leq \mathcal R(\mathcal S) \leq \kappa C$ for a suitable $\kappa>0$ (depending only on $p$). The proposition above $\mathcal R$-boundedness follows from domination.
\begin{cor} \label{domination}
Let  ${\cal S}, {\cal T} \subset B(L^p(\Sigma))$, $1<p<\infty$ and assume that $\cal T$ is $\mathcal R$ bounded and that for every $S \in \cal S$ there exists $T \in \cal T$ such that $|Sf| \leq |Tf|$ pointwise, for every $f \in L^p(\Sigma)$. Then ${\cal S}$ is $\mathcal R$-bounded.
\end{cor}

Let $(A, D(A))$ be a sectorial operator in a Banach space $X$; this means that $\rho (-A) \supset \Sigma_{\pi-\phi}$ for some $\phi <\pi$ and that $\lambda (\lambda+A)^{-1}$ is bounded in $\Sigma_{\pi-\phi}$. The infimum of all such $\phi$ is called the spectral angle of $A$ and denoted by $\phi_A$. Note that $-A$ generates an analytic semigroup if and only if $\phi_A<\pi/2$. The definition of $\mathcal  R$-sectorial operator is similar, substituting boundedness of $\lambda(\lambda+A)^{-1}$ with $\mathcal R$-boundedness in $\Sigma_{\pi-\phi}$. As above one denotes by $\phi^R_A$ the infimum of all $\phi$ for which this happens; since $\mathcal R$-boundedness implies boundedness, we have $\phi_A \le \phi^R_A$.

\medskip

The $\mathcal R$-boundedness of the resolvent characterizes the regularity of the associated inhomogeneous parabolic problem, as we explain now.

An analytic semigroup $(e^{-tA})_{t \ge0}$ on a Banach space $X$ with generator $-A$ has
{\it maximal regularity of type $L^q$} ($1<q<\infty$)
if for each $f\in L^q([0,T];X)$ the function
$t\mapsto u(t)=\int_0^te^{-(t-s)A})f(s)\,ds$ belongs to
$W^{1,q}([0,T];X)\cap L^q([0,T];D(B))$.
This means that the mild solution of the evolution equation
$$u'(t)+Au(t)=f(t), \quad t>0, \qquad u(0)=0,$$
is in fact a strong solution and has the best regularity one can expect.
It is known that this property does not depend on $1<q<\infty$ and $T>0$.
A characterization of maximal regularity is available in UMD Banach spaces, through the $\mathcal  R$-boundedness of the resolvent in a suitable sector $\omega+\Sigma_{\phi}$, with $\omega \in \R$ and $\phi>\pi/2$ or, equivalently, of the scaled semigroup $e^{-(A+\omega')t}$ in a sector around the positive axis. In the case of $L^p$ spaces it can be restated in the following form,  see \cite[Theorem 1.11]{KunstWeis}

\begin{teo}\label{MR} Let $(e^{-tA})_{t \ge0}$ be a bounded analytic semigroup in $L^p(\Sigma)$, $1<p<\infty$,  with generator $-A$. Then $T(\cdot)$ has maximal regularity of type $L^q$  if and only if the set $\{\lambda(\lambda+A)^{-1}, \lambda \in  \Sigma_{\pi/2+\phi} \}$ is $\mathcal R$- bounded for some $\phi>0$. In an equivalent way, if and only if 
there are constants $0<\phi<\pi/2 $, $C>0$ such that for every finite sequence $(\lambda_i) \subset \Sigma_{\pi/2+\phi}$, $(f_i) \subset  L^p$
$$
\left\|\left (\sum_i |\lambda_i (\lambda_i+A)^{-1}f_i|^2\right )^{\frac{1}{2}}\right\|_{L^p(\Sigma)} \le C\left\|\left (\sum_i |f_i|^2\right)^{\frac{1}{2}}\right\|_{L^p(\Sigma)}
$$
or, equivalently, there are constants $0<\phi'<\pi/2 $, $C'>0$ such that  for every finite sequence
$(z_i) \subset \Sigma_{\phi'}$, $(f_i) \subset  L^p$
$$
\left\|\left (\sum_i |e^{-z_i A}f_i|^2\right )^{\frac{1}{2}}\right\|_{L^p(\Sigma)} \le C'\left\|\left (\sum_i |f_i|^2\right)^{\frac{1}{2}}\right\|_{L^p(\Sigma)}.
$$
\end{teo}
\medskip

Finally we state  a version of the operator-valued Mikhlin multiplier theorem in the N-dimensional case, see \cite[Theorem 3.25]{DenkHieberPruss} or \cite[Theorem 4.6]{KunstWeis}.
\begin{teo}   \label{mikhlin}
Let $1<p<\infty$, $M\in C^N(\R^N\setminus \{0\}; B(L^p(\Sigma))$ be such that  the set
$$\left \{|\xi|^{|\alpha|}D^\alpha_\xi M(\xi): \xi\in \R^{N}\setminus\{0\}, \ |\alpha | \leq N \right \}$$
is $\mathcal{R}$-bounded.
Then the operator $T_M={\cal F}^{-1}M {\cal F}$ is bounded in $L^p(\R^N, L^p(\Sigma))$, where $\cal{F}$ denotes the Fourier transform.
\end{teo}

\section{Degenerate operators and similarity transformations }\label{Section Degenerate}
We consider first the 1d operators
$$
L=D_{yy}+\frac{c}{y}D_y-\frac{b}{y^2}, \qquad B=D_{yy}+\frac{c}{y}D_y
$$ on the  half line $\R_+=]0, \infty[$. Note that $B$ (which stands for Bessel) is nothing but $L$ when $b=0$. Often we write $L_y, B_y$ to indicate that they act with respect to the $y$ variable.



The equation $Lu=0$ has solutions $y^{-s_1}$, $y^{-s_2}$ where $s_1,s_2$ are the roots of the indicial equation $f(s)=-s^2+(c-1)s+b=0$ 

\begin{equation} \label{defs}
s_1:=\frac{c-1}{2}-\sqrt{D},
\quad
s_2:=\frac{c-1}{2}+\sqrt{D}
\end{equation}
where
\bigskip
\begin{equation} \label{defD}
D:=
b+\left(\frac{c-1}{2}\right)^2.
\end{equation}

The above numbers are real if and only if $D \ge 0$. When $D<0$ the equation $u-Lu=f$ cannot have positive distributional solutions for certain positive $f$, see \cite{met-soba-spi3}. 
    When $b=0$, then $\sqrt D=|c-1|/2$ and $s_1=0, s_2=c-1$ for $c \ge 1$ and $s_1=c-1, s_2=0$ for $c<1$.

\medskip

%

Next we consider, for $\alpha_1, \alpha_2\in\R$, the  operators
$$\mathcal L=y^{\alpha_1}\Delta_x+y^{\alpha_2} L_y$$ (keeping the assumption  $D\geq 0$ on $ L_y$)
in the space $L^p_m=L^p_m(\R^{N+1}_+)$.

We investigate when these operators can be transformed one into the other by means of  change of variables and multiplications.

For  $k,\beta \in\R$, $\beta\neq -1$ let
\begin{align}\label{Gen Kelvin def}
T_{k,\beta\,}u(x,y)&:=|\beta+1|^{\frac 1 p}y^ku(x,y^{\beta+1}),\quad (x,y)\in\R^{N+1}_+.
\end{align}
Observe that
$$ T_{k,\beta\,}^{-1}=T_{-\frac{k}{\beta+1},-\frac{\beta}{\beta+1}\,}.$$

\begin{prop}\label{Isometry action der} Let $1\leq p\leq \infty$, $k,\beta \in\R$, $\beta\neq -1$. The following properties hold.
	\begin{itemize}
		\item[(i)] For every $m\in\R$,  $T_{k,\beta\,}$ maps isometrically  $L^p_{\tilde m}$ onto $L^p_m$  where 
		$$ \tilde m=\frac{m+kp-\beta}{\beta+1}.$$
		\item[(ii)] For every  $u\in W^{2,1}_{loc}\left(\R^{N+1}_+\right)$ one has
\begin{itemize}
	\item[1.] $y^\alpha T_{k,\beta\,}u=T_{k,\beta\,}(y^{\frac{\alpha}{\beta+1}}u)$, for any $\alpha\in\R$;\medskip
	\item [2.] $D_{x_ix_j}(T_{k,\beta\,}u)=T_{k,\beta} \left(D_{x_ix_j} u\right)$, \quad $D_{x_i}(T_{k,\beta\,}u)=T_{k,\beta}\left(D_{x_i} u\right)$;\medskip
	\item[3.]  $D_y T_{k,\beta\,}u=T_{k,\beta\,}\left(ky^{-\frac 1 {\beta+1}}u+(\beta+1)y^{\frac{\beta}{\beta+1}}D_yu\right)$,
	\\[1ex] $D_{yy} (T_{k,\beta\,} u)=T_{k,\beta\,}\Big((\beta+1)^2y^{\frac{2\beta}{\beta+1}}D_{yy}u+(\beta+1)(2k+\beta)y^{\frac{\beta-1}{\beta+1}}D_y u+k(k-1)y^{-\frac{2}{\beta+1}}u\Big)$.\medskip
	\item[4.] $D_{xy} T_{k,\beta\,}u=T_{k,\beta\,}\left(ky^{-\frac 1 {\beta+1}}D_xu+(\beta+1)y^{\frac{\beta}{\beta+1}}D_{xy}u\right)$
\end{itemize}
	\end{itemize}
\end{prop}{\sc{Proof.}} The proof of (i) follows after observing the Jacobian of $(x,y)\mapsto (x,y^{\beta+1})$ is $|1+\beta|y^{\beta}$. To prove (ii) one can easily observe that any $x$-derivatives commutes with $T_{k,\beta}$. Then we compute  
\begin{align*}
D_y T_{k,\beta\,}u(x,y)=&|\beta+1|^{\frac 1 p}y^{k}\left(k\frac {u(x,y^{\beta+1})} y+(\beta+1)y^\beta D_y u(x,y^{\beta+1})\right)\\[1ex]
=&T_{k,\beta\,}\left(ky^{-\frac 1 {\beta+1}}u+(\beta+1)y^{\frac{\beta}{\beta+1}}D_yu\right)
\end{align*}
and similarly
\begin{align*}
D_{yy} T_{k,\beta\,} u(x,y)=&T_{k,\beta\,}\Big((\beta+1)^2y^{\frac{2\beta}{\beta+1}}D_{yy}u+(\beta+1)(2k+\beta)y^{\frac{\beta-1}{\beta+1}}D_y u+k(k-1)y^{-\frac{2}{\beta+1}}u\Big).
\end{align*}
\qed
\begin{prop}\label{Isometry action}  Let 
	$T_{k,\beta\,}$ be the isometry above defined.   The following properties hold.
	
	\begin{itemize}
		\item[(i)] For every  $u\in W^{2,1}_{loc}\left(\R^{N+1}_+\right)$ one has
		$$
		T_{k,\beta\,}^{-1} \Big(y^{\alpha_1}\Delta_x+y^{\alpha_2} L_y\Big)T_{k,\beta\,}u=\Big(y^{\frac{\alpha_1}{\beta+1}}\Delta_x+(\beta+1)^2y^{\frac{\alpha_2+2\beta}{\beta+1}}\tilde{ L}_y\Big) u
		$$
		where $\tilde { L}$ is the operator defined as in \eqref{defL} with parameters $b,c$ replaced, respectively, by
		\begin{align}
			\label{tilde b}
			\nonumber\tilde b&=\frac{b-k\left(c-1+k\right)}{(\beta+1)^2}
			,\\[1ex]
			\tilde c&=\frac{c+2k+\beta\left(c+1+2k+\beta\right)}{(\beta+1)^2}.
		\end{align}
	\item[(ii)] The discriminant $\tilde D$ and the parameters $\tilde s_{1,2}$ of $\tilde L$ defined as in \eqref{defD}, \eqref{defs}  are given by
	\begin{align}\label{tilde D gamma}
		\tilde D&=\frac{D}{(\beta+1)^2},
	\end{align}
	and
	\begin{align}\label{tilde s gamma}
		\tilde s_{1,2}=\frac{s_{1,2}+k}{\beta+1} \quad  (\beta+1>0), \qquad
		\tilde s_{1,2}&=\frac{s_{2,1}+k}{\beta+1}\quad  (\beta+1<0).
	\end{align}
	\end{itemize}	
\end{prop}
{\sc{Proof.}} Using Proposition \ref{Isometry action der} we can compute
\begin{align*}
	{L}_yT_{k,\beta\,}u(x,y)&=T_{k,\beta\,}\Big[(\beta+1)^2y^{\frac{2\beta}{\beta+1}}D_{yy}u+(\beta+1)(2k+\beta)y^{\frac{\beta-1}{\beta+1}}D_y u+k(k-1)y^{-\frac{2}{\beta+1}}u\\[1ex]
	&\hspace{9ex}+cky^{-\frac 2 {\beta+1}}u+c(\beta+1)y^{\frac{\beta-1}{\beta+1}}D_yu-by^{-\frac{2}{\beta+1}}u \Big]\\[1ex]
	&=T_{k,\beta\,}\Bigg[y^{\frac{2\beta}{\beta+1}}
	\Bigg(
	(\beta+1)^2D_{yy} u+\frac{(\beta+1)\left(2k+\beta+c\right)}{y}D_y u\Bigg.\\[1ex]
	&\Bigg.\hspace{15ex}-\Big(b-k(c+k-1)\Big)\frac u{y^2}
	\Bigg)
	\Bigg]=T_{k,\beta\,}\left(y^{\frac{2\beta}{\beta+1}}\tilde { L_y} u\right)
\end{align*}
which implies
\begin{align*}
	T_{k,\beta\,}^{-1}\left(y^{\alpha_2}{L}_y\right)T_{k,\beta\,}u&=y^{\frac{\alpha_2+2\beta}{\beta+1}}\tilde { L_y} u.
\end{align*}
Similarly one has  $y^{\alpha_1}\Delta_x T_{k,\beta\,}u=T_{k,\beta\,}\left(y^{\frac{\alpha_1}{\beta+1}}\Delta_x u\right)$. Adding the last equalities yields (i). The remaining properties follow  directly from the definitions  \eqref{defs}, \eqref{defD}.
\qed 
\section{Weighted Sobolev spaces} 
Let $p>1$, $m, \alpha_1 \in \R$, $\alpha_2<2$. In order to describe the domain of  the operator $y^{\alpha_{1}}\Delta_x+y^{\alpha_2}B_y$, we collect in this section  the main results concerning anisotropic weighted Sobolev spaces, referring to \cite{MNS-Sobolev} for further details  and all the relative proofs.
We define the Sobolev space
\begin{align*}
W^{2,p}(\alpha_1,\alpha_2,m)&=\left\{u\in W^{2,p}_{loc}(\R^{N+1}_+):\ u,\  y^{\alpha_1} D_{x_ix_j}u,\ y^\frac{\alpha_1}{2} D_{x_i}u,  y^{\alpha_2}D_{yy}u,\ y^{\frac{\alpha_2}{2}}D_{y}u\in L^p_m\right\}
\end{align*}
which is a Banach space equipped with the norm
\begin{align*}
\|u\|_{W^{2,p}(\alpha_1,\alpha_2,m)}=&\|u\|_{L^p_m}+\sum_{i,j=1}^n\|y^{\alpha_1} D_{x_ix_j}u\|_{L^p_m}+\sum_{i=1}^n\|y^{\frac{\alpha_1}2} D_{x_i}u\|_{L^p_m}+\|y^{\alpha_2}D_{yy}u\|_{L^p_m}+\|y^{\frac{\alpha_2}{2}}D_{y}u\|_{L^p_m}.
\end{align*}
Next we add a Neumann boundary condition for $y=0$  in the form $y^{\alpha_2-1}D_yu\in L^p_m$ and set
\begin{align*}
W^{2,p}_{\mathcal N}(\alpha_1,\alpha_2,m)=\{u \in W^{2,p}(\alpha_1,\alpha_2,m):\  y^{\alpha_2-1}D_yu\ \in L^p_m\}
\end{align*}
with the norm
$$
\|u\|_{W^{2,p}_{\mathcal N}(\alpha_1,\alpha_2,m)}=\|u\|_{W^{2,p}(\alpha_1,\alpha_2,m)}+\|y^{\alpha_2-1}D_yu\|_{ L^p_m}.
$$

We  consider also an integral version of the Dirichlet boundary condition, namely  a weighted summability requirement for $y^{-2}u$ and introduce 
$$
	W^{2,p}_{\mathcal R}(\alpha_1, \alpha_2, m)=\{u \in  W^{2,p}(\alpha_1, \alpha_2, m): y^{\alpha_2-2}u \in L^p_m\}
$$
with the norm $$\|u\|_{W^{2,p}_{\mathcal R}(\alpha_1, \alpha_2, m)}=\|u\|_{W^{2,p}(\alpha_1, \alpha_2, m)}+\|y^{\alpha_2-2}u\|_{L^p_m}.$$
The symbol $\mathcal R$ stands for "Rellich", since Rellich inequalities concern with the summability of $y^{-2}u$.

 We consider  only the case  $\alpha_2<2$. Analogous results can be recovered for $\alpha_2>2$ via the similarity transformation of  Lemma \ref{Sobolev eq}.

We have made the choice  not to  include the mixed derivatives in the definition of  $W^{2,p}_{\mathcal{N}}\left(\alpha_1,\alpha_2,m\right)$ to  simplify some arguments. 
However the following result follows from Theorem \ref{complete-Bessel}.

\begin{prop}\label{Sec sob derivata mista} If   $\alpha_2-\alpha_1<2$ and $\alpha_1^{-} <\frac{m+1}p$ then  for every $u \in W^{2,p}_{\mathcal{N}}\left(\alpha_1,\alpha_2,m\right)$  
$$\|y^\frac{\alpha_1+\alpha_2}{2} D_{y}\nabla_x u \|_{ L^p_m} \leq C \|u\|_{W^{2,p}_{\cal N}(\alpha_1, \alpha_2, m)}.$$ 
\end{prop}

\begin{os}\label{Os Sob 1-d}
	With obvious changes we consider also the analogous Sobolev spaces $W^{2,p}(\alpha_2,m)$ and $W^{2,p}_{\cal N}(\alpha_2, m)$ on $\R_+$. 
For example  we have
	$$W^{2,p}_{\mathcal N}(\alpha,m)=\left\{u\in W^{2,p}_{loc}(\R_+):\ u,\    y^{\alpha}D_{yy}u,\ y^{\frac{\alpha}{2}}D_{y}u,\ y^{\alpha-1}D_{y}u\in L^p_m\right\}.$$
	All the results of this section will be valid also in $\R_+$ changing (when it appears) the condition $\alpha_1^{-} <\frac{m+1}p$  to $0<\frac{m+1}p$.
\end{os}

The next proposition shows how these spaces transform under the map of  Section 3.
\begin{prop}
\label{Sobolev eq}
Let $p>1$, $m, \alpha_1,\alpha_2\in \R$ with $\alpha_2< 2$.
 Then one has 
\begin{align*}
W^{2,p}_{\mathcal{N}}(\alpha_1,\alpha_2,m)=T_{0,\beta}\left(W^{2,p}_{\mathcal{N}}(\tilde \alpha_1,\tilde\alpha_2,\tilde m)\right),\qquad \tilde\alpha_1=\frac{\alpha_1}{\beta+1},\quad \tilde\alpha_2=\frac{\alpha_2+2\beta}{\beta+1}.
\end{align*}
In particular, by choosing $\beta=-\frac{\alpha_2}2$ one has 
\begin{align*}
W^{2,p}_{\mathcal{N}}(\alpha_1,\alpha_2,m)=T_{0,-\frac {\alpha_2} 2}\left(W^{2,p}_{\mathcal{N}}(\tilde \alpha,0,\tilde m)\right),\qquad \tilde\alpha=\frac{2\alpha_1}{2-\alpha_2},\quad \tilde m=\frac{m+\frac{\alpha_2} 2}{1-\frac{\alpha_2} 2}.
\end{align*}
\end{prop}

\begin{os}
It is essential to deal with  $W^{2,p}_{\mathcal{N}}(\alpha_1,\alpha_2,m)$: in  general the map $T_{0,\beta}$ does  not  transform  $W^{2,p}(\tilde \alpha_1,\tilde\alpha_2,\tilde m)$ into $W^{2,p}(\alpha_1,\alpha_2,m)$.
\end{os}

The next result clarifies in which sense the condition $y^{\alpha_2-1}D_y u \in L^p_m$ is a Neumann boundary condition.

\begin{prop} \label{neumann} The following assertions hold.
\begin{itemize} 
\item[(i)] If $\frac{m+1}{p} >1-\alpha_2$, then $W^{2,p}_{\mathcal N}(\alpha_1, \alpha_2, m)=W^{2,p}(\alpha_1, \alpha_2, m)$.
\item[(ii)] If $\frac{m+1}{p} <1-\alpha_2$, then $$W^{2,p}_{\mathcal N}(\alpha_1, \alpha_2, m)=\{u \in W^{2,p}(\alpha_1, \alpha_2, m): \lim_{y \to 0}D_yu(x,y)=0\ {\rm for\ a.e.\   x \in \R^N }\}.$$
\end{itemize}
In both cases (i) and (ii), the norm of $W^{2,p}_{\mathcal N}(\alpha_1, \alpha_2, m)$ is equivalent to that of $W^{2,p}(\alpha_1, \alpha_2, m)$.
\end{prop}

We provide  an equivalent description of $W^{2,p}_{\mathcal N}(\alpha_1, \alpha_2, m)$, adapted to the operator $D_{yy}+cy^{-1}D_y$.
\begin{prop}\label{Trace D_yu in W}
		Let   $c\in\R$ and $\frac{m+1}{p}<c+1-\alpha_2$.  Then
		\begin{align*} 
			W^{2,p}_{\mathcal N}(\alpha_1, \alpha_2, m)=&\left\{u \in  W^{2,p}_{loc}(\R^{N+1}_+): u,\ y^{\alpha_1}\Delta_xu\in L^p_m  \right. \\[1ex]			
			&\left.\hspace{10ex} y^{\alpha_2}\left(D_{yy}u+c\frac{D_yu}y\right) \in L^p_m\text{\;\;and\;\;}\lim_{y\to 0}y^c D_yu=0\right\}
		\end{align*}
and the norms $\|u\|_{W^{2,p}_{\mathcal N}(\alpha_1,\alpha_2,m)}$ and $$\|u\|_{L^p_m}+\|y^{\alpha_1}\Delta_x u\|_{L^p_m}+\|y^{\alpha_2}(D_{yy}u+cy^{-1}D_yu)\|_{L^p_m}$$ are equivalent on $W^{2,p}_{\mathcal N}(\alpha_1, \alpha_2, m)$.
		Finally,  when $0<\frac{m+1}p\leq c-1$ then 
				\begin{align*} 
			W^{2,p}_{\mathcal N}(\alpha_1, \alpha_2, m)=&\left\{u \in  W^{2,p}_{loc}(\R^{N+1}_+): u,\ y^{\alpha_1}\Delta_xu,   y^{\alpha_2}\left(D_{yy}u+c\frac{D_yu}y\right) \in L^p_m\right\}.
		\end{align*}
	\end{prop}
The following equivalent description of $W^{2,p}_{\mathcal N}(\alpha_1, \alpha_2, m)$ involves a Dirichlet, rather than Neumann, boundary condition,  in a certain range of parameters.
		\begin{prop}\label{trace u in W op}
			Let   $c\geq 1$ and $\frac{m+1}{p}<c+1-\alpha_2$.  The following properties hold.
			\begin{itemize}
				\item[(i)] If $c>1$ then 
				\begin{align*}
					W^{2,p}_{\mathcal N}(\alpha_1, \alpha_2, m)=&\left\{u \in  W^{2,p}_{loc}(\R^{N+1}_+): u,\ y^{\alpha_1}\Delta_xu\in L^p_m,\right. \\[1ex]			
					&\left.\hspace{11ex}y^{\alpha_2}\left(D_{yy}u+c\frac{D_yu}y\right) \in L^p_m\text{\;and\;}\lim_{y\to 0}y^{c-1} u=0\right\}.
				\end{align*}
			\item[(ii)] If $c=1$ then 
\begin{align*}
	W^{2,p}_{\mathcal N}(\alpha_1, \alpha_2, m)=&\left\{u \in  W^{2,p}_{loc}(\R^{N+1}_+): u,\ y^{\alpha_1}\Delta_xu\in L^p_m,\right. \\[1ex]			
	&\left.\hspace{11ex}y^{\alpha_2}\left(D_{yy}u+c\frac{D_yu}y\right) \in L^p_m\text{\;and\;}\lim_{y\to 0} u(x,y)\in \C\right\}.
\end{align*}
			\end{itemize} 
		\end{prop}
	
The next results  show the density  of smooth functions in $W^{2,p}_{\mathcal N}(\alpha_1,\alpha_2,m)$. Let
\begin{equation} \label{defC}
\mathcal{C}:=\left \{u \in C_c^\infty \left(\R^N\times[0, \infty)\right), \ D_y u(x,y)=0\  {\rm for} \ y \leq \delta\ {\rm  and \ some\ } \delta>0\right \},
\end{equation}
its one dimensional version 
\begin{equation} \label {defD}
\mathcal{D}=\left \{u \in C_c^\infty ([0, \infty)), \ D_y u(y)=0\  {\rm for} \ y \leq \delta\ {\rm  and \ some\ } \delta>0\right \}
\end{equation}
and finally (finite sums below)
 $$C_c^\infty (\R^{N})\otimes\mathcal D=\left\{u(x,y)=\sum_i u_i(x)v_i(y), \  u_i \in C_c^\infty (\R^N), \  v_i \in \cal D \right \}\subset \mathcal C.$$
\begin{teo} \label{core gen}
If  $\frac{m+1}{p}>\alpha_1^-$
then $C_c^\infty (\R^{N})\otimes\mathcal D$ is dense in $W^{2,p}_{\mathcal N}(\alpha_1,\alpha_2,m)$.
\end{teo}

Note that the condition $(m+1)/p>\alpha_1^-$, or $m+1>0$ and $(m+1)/p+\alpha_1>0$, is necessary for the inclusion  $C_c^\infty (\R^{N})\otimes\mathcal D \subset W^{2,p}_{\mathcal N}(\alpha_1,\alpha_2,m)$. 
\medskip

\begin{prop}\label{Hardy Rellich Sob}
	Let   $\frac{m+1}{p}>\alpha_1^-$. The following properties hold for any $u\in W^{2,p}_{\mathcal N}(\alpha_1, \alpha_2, m)$.
	\begin{itemize}
		\item[(i)] If $\frac{m+1}p>1-\frac{\alpha_2}2$ then
		\begin{align*}
			\|y^{\frac{\alpha_2}2-1}u\|_{L^p_m}\leq C \|y^{\frac{\alpha_2}2}D_{y}u\|_{L^p_m}.
		\end{align*}
		\item[(ii)] If  $\alpha_2-\alpha_1<2$ and  $\frac{m+1}p>1-\frac{\alpha_1+\alpha_2}{2}$, $\frac{m+1}p>\alpha_1^-$ then
		\begin{align*}
			\|y^{\frac{\alpha_1+\alpha_2}{2}-1}\nabla_{x}u\|_{L^p_m}\leq C \|y^\frac{\alpha_1+\alpha_2}{2} D_{y}\nabla_x u\|_{L^p_m}.
		\end{align*}
	\end{itemize}
\end{prop}

Finally, we investigate some relationships between $W^{2,p}(\alpha_1, \alpha_2, m)$, $W^{2,p}_{\mathcal R}(\alpha_1, \alpha_2, m)$ and  $W^{2,p}_{\mathcal N}(\alpha_1, \alpha_2, m)$.  

\begin{prop} \label{RN} The following properties hold.
	\begin{itemize}
		\item[(i)]  if $u \in W^{2,p}_{\mathcal R}(\alpha_1, \alpha_2, m)$ then $y^{\alpha_2-1}D_y u \in L^p_m.$
		\item[(ii)] If $\alpha_2-\alpha_1<2$ and $\frac{m+1}{p}>2-\alpha_2$, then 
		$$W^{2,p}_{\mathcal R}(\alpha_1, \alpha_2, m) =  W^{2,p}_{\mathcal N}(\alpha_1, \alpha_2, m)=W^{2,p}(\alpha_1, \alpha_2, m),$$ with equivalence of the corresponding norms. In particular, $C_c^\infty (\R^{N+1}_+)$ is dense in $W^{2,p}_{\mathcal R}(\alpha_1, \alpha_2, m) $.
	\end{itemize}
\end{prop}

We clarifies the action of the multiplication operator $T_{k,0}:u\mapsto y^ku$.  The following lemma  is the companion of Lemma \ref{Sobolev eq} which deals with the transformation $T_{0,\beta}$.
\begin{lem}\label{isometryRN}
	\label{y^k W}
	Let   $\alpha_2-\alpha_1<2$ and  $\frac{m+1}{p}>2-\alpha_2$. For every $k\in\R$
	\begin{align*}
		T_{k,0}:  W^{2,p}_{\mathcal N}(\alpha_1, \alpha_2, m) \to  W^{2,p}_{\mathcal R}(\alpha_1, \alpha_2, m-kp)
	\end{align*}
	is an isomorphism (we shall write $y^k  W^{2,p}_{\mathcal N}(\alpha_1, \alpha_2, m)= W^{2,p}_{\mathcal R}(\alpha_1, \alpha_2, m-kp)$).
\end{lem}
\section{One dimensional degenerate operators}
In this section we summarize  the main   results  proved in \cite{MNS-PerturbedBessel}  for the one dimensional operator  $y^{\alpha} B_y-\mu y^{\alpha}=y^\alpha\left(D_{yy}+\frac{c}{y}D_y \right)-\mu y^\alpha$, $\mu \geq 0$,  in $L^p_m$. To characterize the domain for $\mu>0$,  we denote by $$D(y^\alpha)=\left\{u\in L^p_m:\ y^{\alpha}u\in L^p_m\right\}$$ the domain of the potential $V(y)=y^\alpha$ in $L^p_m$.

\begin{teo}  \label{1d} Let  $\alpha<2$, $c\in\R$ and $1<p<\infty$.    
\begin{itemize}
\item[(i) ]If  $0<\frac{m+1}p<c+1-\alpha$, then  the operator 
		$y^\alpha B$ endowed with domain $W^{2,p}_{\mathcal N}(\alpha,m)$	
		generates a bounded positive analytic semigroup of angle $\pi/2$ on $L^p\left(\R_+,y^mdy\right)$.
\item[(ii)] If $\mu>0$ and $\alpha^-<\frac{m+1}p<c+1-\alpha$ then  the operator 	 
		$y^{\alpha} B-\mu y^{\alpha}$  endowed with domain $W^{2,p}_{\mathcal N}(\alpha,m)\cap D(y^{\alpha})$ generates a bounded analytic semigroup  in $L^p_m$ which has maximal regularity.
\end{itemize}
In both cases the set
	$\mathcal {D}$ defined in \eqref{defD} 
	is a core.
\end{teo}

 We shall use $y^\alpha B^n$, $n$ stands for Neumann,  for $y^\alpha B$ with domain $W^{2,p}_{\mathcal N}(\alpha,m)$ and similarly for $y^\alpha B^n-\mu y^\alpha$. 
  Note that the condition  $\alpha^-<\frac{m+1}p<c+1-\alpha$  is equivalent to   $0<\frac{m+1}p<c+1-\alpha$ and $-\alpha<\frac{m+1}p<c+1-\alpha$. The first guarantees that $y^\alpha B^n$ is a generator in $L^p_{m}$ and the second that $B^n$ is a generator in $L^p_{m+\alpha p}$.

%
In the next proposition we show  that  the multipliers 
\begin{align*}
	\xi \in \R^N \to N_{\lambda}(\xi)&=\lambda (\lambda-y^\alpha B_y+y^\alpha|\xi|^2)^{-1},\\[1ex]
	\xi \in \R^N \to M_{\lambda}(\xi)&=|\xi|^2y^\alpha(\lambda-y^\alpha B_y+y^\alpha|\xi|^2)^{-1}
\end{align*}
satisfy the hypothesis of  Theorem \ref{mikhlin}.
$M_{\lambda}$ is  used in Section \ref{Sect DOm max}  to  characterize the domain of  $\mathcal L=y^\alpha(\Delta_x+B_y)$ whereas  $N_{\lambda}$ to  prove    that $\mathcal L=y^\alpha(\Delta_x+B_y)$ generates an analytic semigroup.

\begin{prop} \label{Rbounded}
	Assume that  $\alpha^-<\frac{m+1}p<c+1-\alpha$.
	Then the families
	\begin{align*}
			\left \{|\xi|^{|\beta|}D^\beta_\xi( M_{\lambda})(\xi): \xi\in \R^{N}\setminus\{0\}, \ |\beta | \leq N ,\lambda \in \C^+\right \},\\[1ex]
			\left \{|\xi|^{|\beta|}D^\beta_\xi(\lambda  N_{\lambda})(\xi): \xi\in \R^{N}\setminus\{0\}, \ |\beta | \leq N ,\lambda \in \C^+\right \}
	\end{align*} 

	are  $\mathcal{R}$-bounded in $L^p_m$.	
\end{prop}

\section{Domain and maximal regularity for $y^\alpha\Delta_x+y^\alpha B^n_y$} \label{Sect DOm max}
Let $c,m\in\R$ and $p>1$. In this section we prove generation results, maximal regularity and   domain characterization for  the  degenerate   operators 
\begin{equation*} 
\mathcal L :=y^\alpha \Delta_x+y^\alpha B^n_y, \quad \alpha<2
\end{equation*}
in $L^p_m\left(\R^{N+1}_+\right)$, where $y^\alpha B^n_y=y^\alpha\left(D_{yy}+\frac{c}{y}D_y \right)$. We start with the $L^2$ theory.

\subsection{The operator $\mathcal L=y^\alpha\Delta_x+y^\alpha B^n_y$ in $L^2_{c-\alpha}$} 
We use the Sobolev spaces of Section 4 and also  $H^{1}_{\alpha,c}:=\{u \in L^2_{c-\alpha} : y^{\frac\alpha 2}\nabla u \in L^2_{c-\alpha}$ equipped with the inner product
\begin{align*}
 \left\langle u, v\right\rangle_{H^1_{\alpha,c}}:= \left\langle u, v\right\rangle_{L^2_{c-\alpha}}+\left\langle y^{\frac\alpha 2}\nabla u,y^{\frac\alpha 2} \nabla v\right\rangle_{L^2_{c-\alpha}}.
 \end{align*}

Let  $\mathcal L$ be  the operator defined on $C_c^\infty(\R^{N+1}_+)$ by 
\begin{align*}
\mathcal L&=y^\alpha\Delta +cy^{\alpha-1}D_y=
y^{-c+\alpha}{\mbox{ div}}\Big(y^{c}\,\nabla u\Big).
\end{align*}
Note that $c=\alpha$ if and only if $\mathcal L$ is formally self-adjoint with respect to the Lebesgue measure. $\mathcal L$ is associated to the  non-negative, symmetric and closed form  in $L^2_{c-\alpha}(\R^{N+1}_+)$ 
\begin{align*}
\mathfrak{a}(u,v)
&:=
\int_{\R^{N+1}_+} \langle y^{\alpha}\nabla u, \nabla \overline{v}\rangle\,y^{c-\alpha} dx\,dy=\int_{\R^{N+1}_+} \langle\nabla u, \nabla \overline{v}\rangle\,y^c dx\,dy=\int_{\R^{N+1}_+} ( -\mathcal Lu)\, \overline{v}\, y^{c-\alpha} dx\,dy,\\[1ex]
D(\mathfrak{a})
&=H^1_{\alpha,c}.
\end{align*}
Accordingly we define the operator with Neumann boundary conditions by
\begin{align} \label{BesselN}
\nonumber D( \mathcal L )&=\{u \in H^1_{\alpha,c}: \exists  f \in L^2_{c-\alpha} \ {\rm such\ that}\  \mathfrak{a}(u,v)=\int_{\R^{N+1}_+} f \overline{v}y^{c-\alpha}\, dz\ {\rm for\ every}\ v\in H^1_{\alpha,c}\},\\  \mathcal Lu&=-f.
\end{align}
By construction  $\mathcal L$ is  a non-positive self-adjoint operator and,  if $u \in D(\mathcal L)$, then  $ u \in H^2_{loc}(\R^{N+1}_+)$ and $ \mathcal Lu=y^\alpha\Delta u+cy^{\alpha-1}D_yu$ by standard arguments. $ \mathcal L$ generates a contractive analytic semigroup $\left\{e^{zy^\alpha \mathcal L}:\ z\in\C_+\right\}$  in $L^2_{c-\alpha}(\R^{N+1}_+)$  and our aim is to characterize its  domain. 

\begin{prop}\label{core Neumnan comp 2}
If   $c+1>|\alpha|$   then 
 the set $C_c^\infty (\R^{N})\otimes \mathcal D$, see \eqref{defD}, 
is a core for $\mathcal L$ in $L^2_{c-\alpha}(\R^{N+1}_+)$.
\end{prop}
{\sc{Proof.}}  We observe, preliminarily, that under the given assumptions  on $\alpha, c$,  the set $C_c^\infty (\R^{N})\otimes \mathcal D$ is contained in $H^1_{\alpha,c}$. Moreover, integrating by parts one sees that any $u\in C_c^\infty (\R^{N})\otimes \mathcal D$ satisfies  \eqref{BesselN} with $\mathcal L u=y^\alpha\Delta_x u+y^\alpha B_y u \in L^2_{c-\alpha}$. This yields  $C_c^\infty (\R^{N})\otimes \mathcal D\subseteq D( \mathcal L )$.

   Since $I-  \mathcal{L}$ is invertible we have to show that $(I-  \mathcal{L})\left(C_c^\infty (\R^{N})\otimes \mathcal D\right)$ is dense in $L^2_{c-\alpha}$ or, equivalently,  that $\left((I- \mathcal{L})\left(C_c^\infty (\R^{N})\otimes \mathcal D\right)\right )^{\perp}=\left\{0\right\}$. Let $v\in L^2_{c-\alpha}(\R^{N+1}_+)$ be such that
\begin{align*}
\int_{\R^{N+1}_+}\left(I- \mathcal{L}\right)f\, \bar{v}\ dx\ y^{c-\alpha} dy
=0, \quad \forall f\in C_c^\infty (\R^{N})\otimes \mathcal D.
\end{align*}
Let us choose $f=a(x)u(y)\in \mathcal{D}$ with  $a\in C_c^\infty (\R^{N})$ and $u\in \mathcal D$.  Taking the  Fourier transform with respect to $x$   we get $\hat f(\xi,y)=\hat a(\xi) u(y)$ and
\begin{align}\label{Core2 eq 1}
\int_{\R^{N+1}_+}\Big[u(y)+ y^{\alpha}|\xi|^2 u(y)-y^\alpha B_y u(y)\Big]\,\hat a(\xi) \ \bar{\hat v}(\xi,y)\ d\xi\, y^{c-\alpha}dy=0.
\end{align}
Fix $\xi_0\in\R^N$, $r>0$ and let $w(\xi)=\frac{1}{|B(\xi_0,r)|}\chi_{B(\xi_0,r)}\in L^2(\R^N)$. Let  $(a_n)_n\in C_c^\infty(\R^{N})$ a sequence of test functions such that $a_n\to \check{w}$ in $L^2(\R^N)$; then  $\hat a_n\to w$ in $L^2(\R^N)$ and writing  \eqref{Core2 eq 1} with $\hat a$ replaced by $\hat a_n$ and letting  $n\to\infty$ we obtain
\begin{align*}
\frac{1}{|B(\xi_0,r)|}\int_{B(\xi_0,r)}d\xi\int_{0}^\infty \Big[u(y)+ y^{\alpha}|\xi|^2 u(y)-y^\alpha B_y u(y)\Big]\,\bar{\hat v}(\xi,y)\ y^{c-\alpha}dy=0.
\end{align*}

Letting $r\to 0$ and using the Lebesgue Differentiation theorem, we have for a.e. $\xi_0\in\R^N$
\begin{align*}
\int_{0}^\infty \Big[u(y)+ y^{\alpha}|\xi_0|^2 u(y)-y^\alpha B_y u(y)\Big]\,\bar{\hat v}(\xi_0,y)\ y^{c-\alpha}dy=0,
\end{align*}
which is valid for every $u\in \mathcal{D}$. Under the given hypotheses  on $c$ and $\alpha$, Theorem \ref{1d}
 implies that $\mathcal{D}$ is a core for the operator $y^\alpha B_y^n-y^\alpha|\xi_0|^2$  in $L^2_{c-\alpha}(\R_+)$. The last equation  then implies $\hat v(\xi_0,\cdot)=0$  for a.e. $\xi_0\in\R^N$ and the proof is complete.


\qed

\begin{teo}\label{Neumnan comp 2}
If  $c+1>|\alpha|$  then 
\begin{align*}
D(\mathcal L)&=W^{2,2}_\mathcal{N}(\alpha,\alpha,c-\alpha)
\end{align*}
\end{teo}
{\sc Proof. } Observe that
\begin{equation*} 
C_c^\infty (\R^{N})\otimes \mathcal D\subset W^{2,2}_\mathcal{N}(\alpha,\alpha,c-\alpha) \cap D(\mathcal L)
\end{equation*}
 and that  it is core for $\mathcal L$ by Proposition \ref{core Neumnan comp 2} and is dense in  $W^{2,2}_\mathcal{N}(\alpha,\alpha,c-\alpha)$ by \ref{core gen}.

We have to show that the graph norm and that of $W^{2,2}_\mathcal{N}(\alpha,\alpha,c-\alpha)$ are equivalent on $C_c^\infty (\R^{N})\otimes \mathcal D$. Since the second is obviously stronger, we have to show the converse.

 We use Proposition \ref{Trace D_yu in W} and  endow  $W^{2,2}_\mathcal{N}(\alpha,\alpha,c-\alpha)$ with the equivalent  norm $$\|u\|_W=\|u\|_{L^2_{c-\alpha}}+\|y^\alpha\Delta_x u\|_{L^2_{c-\alpha}}+\|y^\alpha B_yu\|_{L^2_{c-\alpha}}.$$

Let $u \in  C_c^\infty (\R^{N})\otimes \mathcal D$ and $f= u-\mathcal L u$, so that  $\|u\|_{L^2_{c-\alpha}} \le \|f\|_{L^2_{c-\alpha}}$.     By taking the Fourier transform with respect to $x$ (with co-variable $\xi$) we obtain
\begin{align}\label{Four eq}
(1+|\xi|^2y^\alpha-y^\alpha B^n_y)\hat u(\xi,\cdot)=\hat f(\xi,\cdot), \qquad y^\alpha|\xi|^2 \hat u(\xi, \cdot)=y^\alpha|\xi|^2(1+|\xi|^2y^\alpha-y^\alpha B^n_y)^{-1}\hat f(\xi,\cdot).
\end{align}

This means $y^\alpha\Delta_x u=-{\cal F}^{-1} M(\xi) {\cal F} f$, where ${\cal F}$ denotes the Fourier transform and $M(\xi)=y^\alpha|\xi|^2(1+|\xi|^2y^\alpha-y^\alpha B^n_y)^{-1}$.

The estimate $\|y^\alpha \Delta_x u\|_{L^2_{c-\alpha}} \le C\|f\|_{L^2_{c-\alpha}}$ then follows  from the  boundedness of the multiplier $M$ in $L^2(\R^N; L^2_{c-\alpha}(\R_+))$ which follows from Proposition \ref{Rbounded} and Theorem \ref{mikhlin} and yields  $\|y^\alpha B_y u\|_{L^2_{c-\alpha}} \le C\|f\|_{L^2_{c-\alpha}}$ by difference. 

This gives  the equivalence of the graph norm and of the norm of $W^{2,2}_\mathcal{N}(\alpha,\alpha,c-\alpha)$ on $C_c^\infty (\R^{N})\otimes \mathcal D$ and concludes the proof.
\qed
\subsection{The operator $\mathcal L=y^\alpha\Delta_x+y^\alpha B^n_y$ in $L^p_m$}

In this section we prove  domain characterization and maximal regularity for  the  degenerate   operator 

\begin{equation*} 
\mathcal L =y^\alpha \Delta_x+y^\alpha B^n_y, \quad \alpha<2
\end{equation*}
in $L^p_m$. To avoid any misinterpretation, we often write $ \mathcal{L}_{m,p}$ to emphasize the underlying space on which the operator acts.

We shall use extensively the set $\mathcal{D}$ defined in \eqref{defD}. In particular $\mathcal L$ is well defined on $C_c^\infty (\R^{N})\otimes\mathcal D$ when $(m+1)/p >\alpha^-$. 

\begin{lem}\label{lemma multipl}
Let $\alpha^-<\frac{m+1}p<c+1-\alpha$.
Then for any $\lambda\in \C^+$ the operators  
$$(\lambda- \mathcal{L}_{c-\alpha,2})^{-1},\quad y^\alpha \Delta_x (\lambda- \mathcal{L}_{c-\alpha,2})^{-1},\quad y^\alpha B^n_y (\lambda- \mathcal{L}_{c-\alpha,2})^{-1}$$ initially defined on $L^p_m\cap L^2_{c-\alpha}$ by Theorem \ref{Neumnan comp 2}, extend to bounded operators on $L^p_m$ which we denote respectively  by $\mathcal{R}(\lambda)$, $y^\alpha \Delta_x \mathcal{R}(\lambda)$, $y^\alpha B^n_y \mathcal{R}(\lambda)$. Moreover  the family  $\left\{\lambda \mathcal{R}(\lambda):\lambda\in\C^+\right\}$ is $\mathcal{R}$-bounded on $L^p_m$.
\end{lem}
{\sc Proof.}
Let   $u \in  C_c^\infty (\R^{N})\otimes  \mathcal{D}$ and $f=\lambda u- \mathcal Lu$.   By taking the Fourier transform with respect to $x$ we obtain
$$
(\lambda+|\xi|^2y^\alpha -y^\alpha B^n_y)\hat u(\xi,\cdot)=\hat f(\xi,\cdot), \qquad  \hat u(\xi, \cdot)=(\lambda -y^\alpha B^n_y+|\xi|^2y^\alpha)^{-1}\hat f(\xi,\cdot).
$$
This means $u={\cal F}^{-1} N_\lambda(\xi) {\cal F} f$, where 
$$N_\lambda(\xi)=(\lambda -y^\alpha B^n_y+|\xi|^2y^\alpha)^{-1}.$$
Since $ C_c^\infty (\R^{N})\otimes  \mathcal{D}$ is a core for $\mathcal{L}_{c-\alpha,2}$ we have proved the equality
$$(\lambda-\mathcal{L}_{c-\alpha,2})^{-1}={\cal F}^{-1} N_\lambda(\xi) {\cal F}.$$
 Proposition \ref{Rbounded} and  Theorem \ref{mikhlin} yield the boundedness  of the Fourier multiplier  $N_\lambda$  in the space $L^p\left(\R^N, L^p_m(\R_+)\right)=L^p_m$ and the existence of a bounded operator $\mathcal{R}(\lambda)\in L^p_m$ which extends $(\lambda- \mathcal{L}_{c-\alpha,2})^{-1}$. Furthermore \cite[Theorem 4.3.9]{Pruss-Simonett} and the $\mathcal R$-boundedness with respect to $\lambda$ of $N_\lambda(\xi)$ and its $\xi$-derivatives, see again Proposition \ref{Rbounded},  imply that  the family $\left\{\lambda \mathcal{R}(\lambda):\lambda\in\C^+\right\}$ is $\mathcal{R}$-bounded.

 The proof for $y^\alpha \Delta_x \mathcal{R}(\lambda)$ is similar. As before we show that, see \eqref{Four eq} in  Theorem \ref{Neumnan comp 2}, 
$$y^\alpha \Delta_x (\lambda- \mathcal{L}_{c-\alpha,2})^{-1}=-{\cal F}^{-1} M_\lambda(\xi) {\cal F}$$ where 
$M_\lambda(\xi)=y^\alpha|\xi|^2(\lambda+|\xi|^2y^\alpha-y^\alpha B^n_y)^{-1}$, and use Proposition \ref{Rbounded}
 for the  boundedness of the multiplier $M_\lambda$ in $L^p(\R^N; L^p_{m}(\R_+))$. 

The boundedness  of $y^\alpha B^n_y\mathcal{R}(\lambda)$ follows then by difference, since $y^\alpha \Delta_x\mathcal{R}(\lambda)+y^\alpha B^n_y\mathcal{R}(\lambda)=\lambda \mathcal{R}(\lambda)-I$.
\qed

\begin{prop} \label{generazione} If  $ \alpha^-<\frac{m+1}p<c+1-\alpha$, an extension $\mathcal L_{m,p}$ of the operator $\mathcal L$, initially defined on $C_c^\infty (\R^{N})\otimes  \mathcal{D}$,  generates a bounded analytic semigroup  in $L^p_m(\R_+^{N+1})$ which has maximal regularity and it is consistent with the semigroup generated by $\mathcal L_{c-\alpha,2}$ in $L^2_{c-\alpha}(\R^{N+1}_+)$. 
\end{prop}
{\sc Proof.}
Let us consider the $\mathcal{R}$-bounded  family of operators $\left\{\lambda \mathcal{R}(\lambda):\lambda\in\C^+\right\}$ defined by Lemma \ref{lemma multipl}. In particular it  satisfies
 \begin{align*}
 \|\lambda \mathcal{R}(\lambda)\|_{\mathcal{B}(L^p_m(\R_+^{N+1}))}\leq C,\qquad\forall \lambda\in\C^+.
 \end{align*}
 By construction  $\mathcal{R}({\lambda})$ coincides with  $(\lambda-\mathcal{L}_{c-\alpha,2})^{-1}$ when restricted to $L^p_m \cap L^2_{c-\alpha}$.  Hence, by density,  the family $\left\{\mathcal{R}(\lambda):\lambda\in\C^+\right\}$  satisfies  the resolvent equation
\begin{align*}
\mathcal{R}(\lambda)-\mathcal{R}(\mu)=(\mu-\lambda)\mathcal{R}(\lambda)\mathcal{R}(\mu),\quad \forall \lambda,\mu\in\C^+
\end{align*}
in $L^p_m$ and therefore it is a pseudoresolvent, see  \cite[Section 4.a]{engel-nagel}. Furthermore $\mbox{rg}(\mathcal{R}(\lambda))$  is dense in $L^p_m$ for every $\lambda \in \C^+$,  since it contains $ C_c^\infty (\R^{N})\otimes  \mathcal{D}$. 

Let us prove that $\mathcal{R}(\lambda)$ is injective for every $\lambda\in \C^+$. Let $f\in L^p_m$ s.t. $\mathcal{R}(\lambda)f=0$ for some $\lambda \in \C^+$. Since   $\mbox{Ker}(\mathcal{R}(\lambda))=\mbox{Ker}(\mathcal{R}(\mu))$ for any $\lambda,\mu\in\C^+$, see \cite[Lemma 4.5]{engel-nagel},  we have  $\mathcal{R}(\lambda)f=0$ for every $\lambda>0$. Given $\epsilon>0$,  let us choose $g\in L^p_m\cap L^2_{c-\alpha}$ s.t. $\|f-g\|_{L^p_m}<\epsilon$. Then  
\begin{align*}
\lambda R(\lambda)g=\lambda R(\lambda)(g-f),\qquad \|\lambda R(\lambda)g\|_{L^p_m}\leq C\epsilon, \qquad \forall\lambda>0.
\end{align*}
Since $\lambda R(\lambda)g=\lambda(\lambda-\mathcal{L}_{c-\alpha,2})^{-1}g\to g$ as $\lambda\to \infty$ we may suppose, up to  a  subsequence, that $\lambda R(\lambda)g\to g$ a.e.. Then Fatou's Lemma yields
\begin{align*}
\|g\|_{L^p_m}\leq \liminf_{\lambda\to\infty} \|\lambda \mathcal{R}(\lambda)g\|_{L^p_m}\leq C\epsilon
\end{align*}
which implies  $\|f\|_{L^p_m}\leq \|f-g\|_{L^p_m}+\|g\|_{L^p_m}\leq (1+C)\epsilon$, hence  $f=0$ which proves the injectivity of $\mathcal{R}(\lambda)$. 

At this point, \cite[Proposition 4.6]{engel-nagel} yields the existence of a densely defined closed operator $\mathcal L_{m,p}$  such that $\C^+\subseteq \rho(\mathcal L_{m,p})$ and $\mathcal R(\lambda)=(\lambda-\mathcal L_{m,p})^{-1}$ for any $\lambda\in\C^+$. By construction, $(\mathcal L_{m,p};D(\mathcal L_{m,p}))$ extends $\left(\mathcal L, C_c^\infty (\R^{N})\otimes  \mathcal{D}\right)$ and one has
 \begin{align*}
 \|\lambda\left(\lambda-\mathcal L_{m,p}\right)^{-1}\|_{\mathcal B\left( L^p_m \right)}\leq C,\qquad \lambda\in C^+.
 \end{align*}
Then from standard results on semigroup theory, see for example \cite[Section AII, Theorem 1.14]{nagel}, $(\mathcal L_{m,p},D(\mathcal L_{m,p}))$ generates a bounded analytic semigroup  $\left(e^{z\mathcal{L}_{m,p}}\right)_{z\in \Sigma_{\theta}}$ for some $\theta>0$,  in $L^p_m$.

The  maximal regularity of the semigroup follows, using Theorem \ref{MR}, from  the $\mathcal{R}$-boundedness of  the resolvent family $\{\lambda\left(\lambda-\mathcal L_{m,p}\right)^{-1},\  \lambda \in \C^+ \}$. Finally, the semigroup is consistent with that in $L^2_{c-\alpha}$ since the resolvents are consistent.
\qed

 Finally we characterize the domain of $\mathcal L_{m,p}$.
\begin{teo} \label{generazione1} 
If $\alpha^-<\frac{m+1}p<c+1-\alpha$, then
\begin{align*}
D(\mathcal L_{m,p})&=W^{2,p}_{\mathcal{N}}(\alpha,\alpha,m)
\end{align*}
 and in particular $C_c^\infty (\R^{N})\otimes \mathcal D$ is a core for $\mathcal L_{m,p}$.
\end{teo}
{\sc Proof.} With the notation of the above proposition, $D(\mathcal L_{m,p})=R(1)\left( L^p_m\right )$.  Let $u=R(1)f= (I- \mathcal{L}_{c-\alpha,2})^{-1}f$ with $f \in L^2_{c-\alpha}\cap L^p_m$. Then 
Lemma \ref{lemma multipl} yields
\begin{align}\label{eq spez}
\|y^\alpha \Delta_x u\|_{L^p_{m}} +\|y^\alpha B_y u\|_{L^p_{m}}\leq C\left(\| \mathcal L u\|_{L^p_{m}}+\| u\|_{L^p_{m}}\right).
\end{align}

Using Theorem \ref{1d} and  Theorem \ref{Neumnan comp 2}, we deduce that  $u(x,\cdot)\in D(y^\alpha B^n_{c-\alpha,2})$ for a.e. $x\in \R^n$. Moreover, $u(x,\cdot)$, $\ y^\alpha B_y u(x, \cdot)\in L^p_m(\R_+)$,  for a.e. $x\in \R^n$. 

Let us show that  $u(x,\cdot)\in D(y^\alpha B^n_{m,p})$. In fact,   setting $f:=u(x,\cdot)-B_y u(x,\cdot)\in L^p_m(\R_+)\cap L^2_{c-\alpha}(\R_+)$ we have $u=\left(I-y^\alpha B^n\right)^{-1}f\in D(y^\alpha B^n_{m,p})\cap D(y^\alpha B^n_{c-\alpha,2})$ by the consistency of the resolvent $\left(I-y^\alpha B^n\right)^{-1}$ in $L^p_m(\R_+)$ and  in $L^2_{c-\alpha}(\R_+)$ .  

Theorem \ref{1d} then implies
\begin{align*}
\|y^{\alpha} D_{yy}u\|_{L^p_{m}(\R_+)}+\|y^{\frac{\alpha}2} D_{y}u\|_{L^p_{m}(\R_+)}+ \|y^{\alpha_2} D_{y}u\|_{L^p_{m}(\R_+)}\leq C\|u-y^\alpha  B_y u\|_{L^p_{m}(\R_+)}.
\end{align*}
Then, raising to the power $p$,  integrating over $\R^N$ and using Lemma \ref{lemma multipl} for the last inequality 
\begin{align} \label{eqspez1}
\|y^{\alpha}D_{yy}u\|_{L^p_{m}}+\|y^{\frac{\alpha}{2}}D_{y}u\|_{L^p_{m}}+\|y^{\alpha-1}D_yu\|_{L^p_{m}}
\leq C \left( \|u\|_{L^p_{m}}+\|y^{\alpha}B_y u\|_{L^p_{m}}\right )
\leq C\left( \|u\|_{L^p_{m}}+\|\mathcal {L} u\|_{L^p_{m}}\right ).
\end{align}

By the density of  $L^2_{c-\alpha}\cap L^p_m$ in $L^p_m$, \eqref{eq spez}, \eqref{eqspez1} hold for every $u \in D(\mathcal L_{m,p})$ and this last is contained in  $W^{2,p}_\mathcal{N}(\alpha,\alpha,m)$, by  \ref{Trace D_yu in W}.

Moreover, since the graph norm is clearly weaker than the norm of $W^{2,p}_\mathcal{N}(\alpha,\alpha,m)$,  \eqref{eq spez}, \eqref{eqspez1} again show that they are equivalent on 
$ D(\mathcal L_{m,p})$, in particular on $C_c^\infty (\R^{N})\otimes \mathcal D$ which is dense in $W^{2,p}_\mathcal{N}(\alpha,\alpha,m)$, by \ref{core gen}. 

Therefore $ D(\mathcal L_{m,p})=W^{2,p}_\mathcal{N}(\alpha,\alpha,m)$ and in particular $C_c^\infty (\R^{N})\otimes \mathcal D$ is a core.
\qed

\begin{cor} \label{omogeneo}
Under the hypotheses of Theorem \ref{generazione1} we have for every $u \in W^{2,p}_{\mathcal{N}}(\alpha,\alpha,m)$
$$
\|y^\alpha D_{x_i x_j} u\|_{L^p_{m}} +\|y^\alpha D_{yy} u\|_{L^p_{m}}+\|y^{\alpha-1} D_{y} u\|_{L^p_{m}}\leq C\| \mathcal Lu\|_{L^p_{m}}.
$$
\end{cor}
{\sc Proof.} By Theorem \ref{generazione1} the above inequality holds if $\|  u\|_{L^p_{m}(\R^{N+1}_+)}$ is added to the right hand side.
Applying it  to $u_\lambda (x,y)=u(\lambda x, \lambda y)$, $\lambda >0$ we obtain
\begin{align*}
\|y^\alpha D_{x_i x_j} u\|_{L^p_{m}} +\|y^\alpha D_{yy} u\|_{L^p_{m}}+\|y^{\alpha-1} D_{y} u\|_{L^p_{m}}\leq C\left(\| \mathcal L u\|_{L^p_{m}}+\lambda^{\alpha-2}\| u\|_{L^p_{m}}\right)
\end{align*}
and the proof follows letting $\lambda \to \infty$.
\qed
\subsection{Mixed  derivatives}
By using classical covering results, Rellich inequalities and Theorem \ref{generazione1}, we obtain   $L^p$ estimates for the mixed second order derivatives.

\begin{teo}   \label{mixed-derivatives}
Let  $\alpha^{-}<\frac{m+1}{p}<c+1-\alpha$. Then there exists $C>0$ such that for every  $u\in D(\mathcal {L}_{m,p})$
$$\|y^\alpha D_{y}\nabla_x u\|_{L^p_m}\leq C\|\mathcal Lu\|_{L^p_m}.$$
\end{teo}

We need a Rellich type inequality for smooth functions vanishing near $\{y=0 \}$.

\begin{lem} \label{Rellich}   
Let $\alpha^{-}<\frac{m+1}{p}<c+1-\alpha$. Assume, in addition, $\alpha\neq 1-\frac{m+1}{p}$, $\alpha\neq 2-\frac{m+1}{p}$. 
Then there exists a positive constant $C$ such that for $u\in  C_c^\infty (\R^{N} \times ]0, \infty[)$ we have
$$\left\|y^{\alpha-2} u\right\|_{L^p_m}\leq C\|\mathcal Lu\|_{L^p_m}.$$
\end{lem}
{\sc Proof.} 
Let  $u\in  C_c^\infty (\R^{N} \times ]0, \infty[)$. Let $\alpha\neq 1-\frac{m+1}{p}$, $\alpha\neq 2-\frac{m+1}{p}$.  Then by  \cite[Proposition 3.10]{met-soba-spi-Rellich} (see also \cite{met-negro-soba-spina})
$$\int_{\R^+}\left|y^{\alpha-2} u\right|^p\, y^m dy\leq C\int_{\R^+}|y^\alpha D_{yy} u|^p\, y^m dy.$$ Integrating the previous inequality over $\R^N$ and  using   Corollary  \ref{omogeneo} we get
$$\left\|y^{\alpha-2} u\right\|_{L^p_m}\leq C \|y^\alpha D_{yy}u\|_{L^p_m}\leq C\|\mathcal Lu\|_{L^p_m}.$$
\qed

We first prove mixed derivatives estimates for functions with support far away from $\{y=0\}$. 
\begin{lem} \label{boundaryAss}
Let $\alpha^{-}<\frac{m+1}{p}<c+1-\alpha$. Assume, in addition, $\alpha\neq 1-\frac{m+1}{p}$, $\alpha\neq 2-\frac{m+1}{p}$.  Then for every $u\in  C_c^\infty (\R^{N} \times ]0, \infty[)$ 
$$\|y^\alpha D_{y}\nabla_x u\|_{L^p_m}\leq C\|\mathcal Lu\|_{L^p_m}.$$
\end{lem}
{\sc Proof.} 
For every $n\in\Z$ let 
\begin{align*}
I_n=[2^n, 2^{n+1}[, \qquad J_n=[2^{n-1}, 2^{n+2}[
\end{align*}
We fix $\vartheta\in C_c^{\infty }(\R)$ such that $0\leq \vartheta \leq 1$, $\vartheta(y)=1$ for $y\in [1,2]$ and $\vartheta(y)=0$ for $y\notin [\frac 12,4]$ and  set $\vartheta_n(y)=\vartheta \left(\frac{y}{\rho_n}\right)$, where $\rho_n=2^{n}$.

We apply the classical $L^p$  estimates for elliptic operators with constant coefficients
to the function $\vartheta_n u$ and obtain
$$\|\rho_n^\alpha D_{y}\nabla_x(\vartheta_n u)\|_{L^p(\R^{N+1}_+)}\leq C\|\rho_n^\alpha D_{yy}(\vartheta_n u)+\rho_n^\alpha\Delta_x(\vartheta_n u)\|_{L^p(\R^{N+1}_+)}.$$
Then we get
\begin{align*}
\|\rho_n^\alpha D_{y}\nabla_x u\|_{L^p(\R^N\times I_n)}&\quad\leq C\Big(\|\rho_n^\alpha D_{yy}u+\rho_n^\alpha\Delta_x u\|_{L^p(\R^N\times J_n)}+\frac{1}{\rho_n}\|\rho_n^\alpha D_y u\|_{L^p(\R^N\times J_n)}\\&\quad+\frac{1}{\rho_n^2}\|\rho_n^\alpha u\|_{L^p(\R^N\times J_n}\Big)\\[1ex]
\end{align*}
Since  $\frac{\rho_n}{2}\leq y\leq 4\rho_n$ if $y\in J_n$ then we get
\begin{align*}
\|y^{\alpha+\frac{m}{p}} D_{y}\nabla_x u\|_{L^p(\R^N\times I_n)}   
&\leq C\Big(\|y^{\alpha+\frac{m}{p}} D_{yy}u+y^{\alpha+\frac{m}{p}}\Delta_x u\|_{L^p(\R^N\times J_n)}\\
&+\|y^{\alpha-1+\frac{m}{p}} D_y u\|_{L^p(\R^N\times J_n)}+\|y^{\alpha-2+\frac{m}{p}} u\|_{L^p(\R^N\times J_n)}\Big).
\end{align*}
Summing over $n$, since at most three among the  intervals $J_n$  overlap, it follows  that
\begin{align*}
\|y^\alpha D_{y}\nabla_x u\|_{L^p_m}&\leq C\left(\|\mathcal L_{m,p} u\|_{L^p_m}+\|y^{\alpha-1}D_y u\|_{L^p_m}+\left\|y^{\alpha-2}u\right\|_{L^p_m}\right).
\end{align*}
Using Corollary \ref{omogeneo} and  Lemma \ref{Rellich} we conclude the proof.
\qed

Next we remove the assumption  on the supports and work in $C_c^\infty(\R^N)\otimes \mathcal{D}$ which is a core for $\mathcal L_{m,p}$.
\begin{lem}   \label{boundarAssum2}
Let $\alpha^{-}<\frac{m+1}{p}<c+1-\alpha$ and assume also that  $\alpha\neq 1-\frac{m+1}{p}$, $\alpha\neq 2-\frac{m+1}{p}$. Then
$$\|y^\alpha D_{y}\nabla_xu\|_{L^p_m}\leq C\|\mathcal{L} u\|_{L^p_m}$$ for every $u\in C_c^\infty(\R^N)\otimes \mathcal D$.
\end{lem}
{\sc Proof.} 
Given  $u\in  C_c^\infty(\R^N)\otimes \mathcal D$,  let $v(x,y)=u(x,\lambda y)$. Then $v \in  C_c^\infty(\R^N)\otimes \mathcal D$ and  $u(x,0)=v(x,0)$. It follows that $w=u-v \in  C_c^\infty (\R^{N} \times ]0, \infty[)$. Moreover
$$\|y^\alpha D_{y} \nabla_xv\|_{L^p_m}=
\lambda ^{1-\alpha-\frac{m+1}{p}} \|y^\alpha D_{y} \nabla_x u\|_{L^p_m}$$ and, by Corollary \ref{omogeneo}, 
$$\|\mathcal L v\|_{L^p_m}\leq \lambda^{-\alpha-\frac{m+1}{p}}\|y^\alpha\Delta_x u\|_{L^p_m}+\lambda^{2-\alpha-\frac{m+1}{p}}\|y^\alpha B^n_y u\|_{L^p_m}\leq C(\lambda) \|\mathcal L u\|_{L^p_m}.$$
 Hence by applying Lemma \ref{boundaryAss} to $w$, we have 
\begin{align*}
&\|y^\alpha D_{y} \nabla_x u\|_{L^p_m}\leq C\left(\|y^\alpha D_{y} \nabla_x w\|_{L^p_m}+\|y^\alpha D_{y} \nabla_x v\|_{L^p_m}\right)\leq C\left(\|\mathcal L w\|_{L^p_m}+\|y^\alpha D_{y} \nabla_x v\|_{L^p_m}\right)\\&\leq C\left(\|\mathcal L u\|_{L^p_m}+\|\mathcal L v\|_{L^p_m}+\|y^\alpha D_{y} \nabla_x v\|_{L^p_m}\right)\leq C'(\lambda)\|\mathcal{L}u\|_{L^p_m}+ C \lambda ^{1-\alpha-\frac{m+1}{p}} \|y^\alpha D_{y} \nabla_x u\|_{L^p_m}.
\end{align*}
Choosing $\lambda$ large enough or small enough accordingly to  $1-\alpha-\frac{m+1}{p}>0$ or $1-\alpha-\frac{m+1}{p}<0$ we conclude the proof. \qed

{\sc Proof.} (Theorem \ref{mixed-derivatives}).
Since $ C_c^\infty(\R^N)\otimes \mathcal D$ is a core for $\mathcal L_{m,p}$, by Lemma \ref{boundarAssum2} the claim holds for $\alpha\neq 1-\frac{m+1}{p}$, $\alpha\neq 2-\frac{m+1}{p}$. 

Suppose now   $p=\frac{m+1}{1-\alpha}$ (in particular $\alpha<1$ and $m+\alpha>0$). Observe that, by the previous part of the proof, the operator $y^\alpha D_{y}\nabla_x(I-\mathcal L_{m,q})^{-1}$ is bounded in $L^{q}$ for $q<\frac{m+1}{1-\alpha}$ and for $q> \frac{m+1}{1-\alpha}$, $q$ close to $p$ to satisfy the condition $\alpha^-<(m+1)/q<c+1-\alpha$.  The Riesz-Thorin interpolation theorem then yields the boundedness of $y^\alpha D_{y}\nabla_x(I-\mathcal L_{m,p})^{-1}$ also for $p= \frac{m+1}{1-\alpha}$. By arguing similarly for $p=\frac{m+1}{2-\alpha}$ we conclude the proof.  
 \qed


\section{The operator $y^{\alpha_1}\Delta_x+y^{\alpha_2} B^n_y$}\label{Sec comp alpha N}
In this section we consider for $\alpha_1\in\R$, $\alpha_2<2$  the  operator  $$\mathcal L^{\alpha_1,\alpha_2}_{m,p}=y^{\alpha_1}\Delta_x+y^{\alpha_2} B^n_y$$ in the space $L^p_m$. 
The generation and domain properties for $\mathcal L^{\alpha_1,\alpha_2}$ are deduced from the case $\alpha_1=\alpha_2$ by using the isometry  
\begin{align*}\label{Gen Kelvin def}
T_{k,\beta\,}u(x,y)&:=|\beta+1|^{\frac 1 p}y^ku(x,y^{\beta+1}),\quad (x,y)\in\R^{N+1}_+
\end{align*}
introduced in Section \ref{Section Degenerate}.

\begin{teo} \label{complete-Bessel} 
Let  $\alpha_2-\alpha_1<2$ and
$$\alpha_1^{-} <\frac{m+1}p<c+1-\alpha_2.$$   
  Then $\mathcal L^{\alpha_1,\alpha_2}$ with domain $D(\mathcal L^{\alpha_1,\alpha_2}_{m,p})= W^{2,p}_{\mathcal{N}}\left(\alpha_1,\alpha_2,m\right)$ generates a bounded analytic semigroup  in $L^p_m$ which has maximal regularity.
Moreover for every  $u\in W^{2,p}_{\mathcal{N}}\left(\alpha_1,\alpha_2,m\right)$
$$\|y^\frac{\alpha_1+\alpha_2}{2} D_{y}\nabla_x u\|_{L^p_m}\leq C\|\mathcal Lu\|_{L^p_m}.$$
\end{teo}
{\sc Proof.} We use the isometry 
$$T_{0,\frac{\alpha_1-\alpha_2}{2}}:L^p_{\tilde{m}}\to L^p_m,\qquad \tilde m=\frac{2m-\alpha_1+\alpha_2}{\alpha_1-\alpha_2+2}$$
which, according to Proposition \ref{Isometry action},   transforms  $\mathcal L^{\alpha_1,\alpha_2}$ into  
$$
T_{0,\frac{\alpha_1-\alpha_2}{2}}^{-1} \mathcal L^{\alpha_1,\alpha_2} T_{0,\frac{\alpha_1-\alpha_2}{2}}=y^\alpha\Delta_x+\left(\frac{\alpha_1-\alpha_2+2}{2}\right)^2y^\alpha \tilde B^n_y
$$ where $$\alpha=\frac{2\alpha_1}{\alpha_1-\alpha_2+2}, \quad \tilde B^n_y=D_{yy}+\frac{\tilde c}{y}D_y ,\quad
 \tilde c=\frac{4c+(\alpha_1-\alpha_2)(2c+2+\alpha_1-\alpha_2)}{(\alpha_1-\alpha_2+2)^2}.$$
Observe that $\alpha<2$ by assumption as well as $\alpha^- < \frac{\tilde m +1}{p} < \tilde c+1-\alpha$.
Generation properties and maximal regularity  for $\mathcal L^{\alpha_1\alpha_2}_{m,p}$ in $L^p_m$ are then  immediate consequence of the same properties of $$y^\alpha\Delta_x+\left(\frac{\alpha_1-\alpha_2+2}{2}\right)^2y^\alpha\tilde B^n_y$$ in $L^p_{\tilde m}$ proved in  Proposition \ref{generazione} and Theorem \ref{generazione1}. 
Concerning the domain, we have
$$D(\mathcal L^{\alpha_1,\alpha_2}_{m,p})=T_{0,\frac{\alpha_1-\alpha_2}{2}}\left( W^{2,p}_{\mathcal{N}}\left (\alpha, \alpha,  \tilde m \right)\right)$$
which, by \ref{Sobolev eq}, coincides with $W^{2,p}_{\mathcal{N}}\left(\alpha_1,\alpha_2,m\right)$.
The estimates for the  mixed derivatives follow from the equality
 $$y^{\frac{\alpha_1+\alpha_2}{2}}D_{xy}u=\frac{2+\alpha_1-\alpha_2}{2}T_{0,\frac{\alpha_1-\alpha_2}{2}}\left(y^{\frac{2\alpha_1}{\alpha_1-\alpha_2+2}}D_{xy}\tilde u\right).$$
and Theorem \ref{mixed-derivatives}.
\qed

\begin{os}
 The operator $y^\alpha\Delta_x+ay^\alpha B^n_y$, $a>0$,  has the same domain and properties of $y^\alpha\Delta_x+y^\alpha B^n_y$. This follows by using the map $Tu(x,y)=u(x, a^{-\frac 1 2}y)$ since $T^{-1}\left(y^\alpha\Delta_x+ay^\alpha B^n_y\right)T =a^{\frac\alpha 2}\left(y^\alpha\Delta_x+y^\alpha B^n_y\right)$. We  used this in the above proof.
\end{os}



\section{Degenerate operators with Dirichlet boundary conditions}\label{sec diric}
In this section we add a potential term to $B$ and study  the  operator  
$$
\mathcal{L}=\mathcal{L}^{\alpha_1, \alpha_2}=y^{\alpha_1}\Delta_x+y^{\alpha_2}L_y=y^{\alpha_1}\Delta_x+y^{\alpha_2}\left(D_{yy}+\frac{c}{y}D_y-\frac{b}{y^2}\right), \quad \alpha_2<2
$$
in $L^p_m$, under Dirichlet boundary conditions, in the sense specified below.

We recall that the equation $L_yu=0$ has solutions $y^{-s_1}$, $y^{-s_2}$ where $s_1,s_2$ are the roots of the indicial equation $f(s)=-s^2+(c-1)s+b=0$ given by

\begin{equation*} 
s_1:=\frac{c-1}{2}-\sqrt{D},
\quad
s_2:=\frac{c-1}{2}+\sqrt{D}
\end{equation*}
where
\bigskip
\begin{equation*} 
D:=
b+\left(\frac{c-1}{2}\right)^2
\end{equation*}
is supposed to be nonnegative. 
 When $b=0$, then $\sqrt D=|c-1|/2$ and $s_1=0, s_2=c-1$ for $c \ge 1$ and $s_1=c-1, s_2=0$ for $c<1$.

\begin{os} 
	All the results of this section will be valid, with obvious changes,  also  in  $\R_+$ for the $1$d operators $y^{\alpha_2} L_y$ changing (when it appears in the various conditions on the parameters) $\alpha_1^{-}$ to $0$ (see also Remark  \ref{Os Sob 1-d}).  We  also refer to   \cite{met-calv-negro-spina, MNS-Sharp, MNS-Grad, MNS-Max-Reg,  Negro-Spina-Asympt} for the analogous results concerning the $Nd$ version of $L_y$.
\end{os}

\medskip
A  multiplication operator transforms $\mathcal L$  into an operator of the form $y^{\alpha_1}\Delta_x+y^{\alpha_2}B^n_y$  and allows  to transfer the results of the previous sections to this  situation. 
Indeed,  we  use the map defined in Section \ref{Section Degenerate}
\begin{align}\label{Gen Kelvin def}
T_{k,0\,}u(x,y)&:=y^ku(x,y),\quad (x,y)\in\R^{N+1}_+
\end{align}
for a suitable choice of $k$ and with $\beta=0$.
We recall that  $T_{k,0\,}$ maps isometrically  $L^p_{\tilde m}$ onto $L^p_m$  where 
$ \tilde m=m+kp$ and for every
 $u\in W^{2,1}_{loc}\left(\R^{N+1}_+\right)$ one has
$$
T_{k,0\,}^{-1} \Big(y^{\alpha_1}\Delta_x+y^{\alpha_2} L_y\Big)T_{k,0\,}u=\Big(y^{\alpha_1}\Delta_x+y^{\alpha_2}\tilde{L}_y\Big) u
$$
where $\tilde {L}$ is the operator defined as above with parameters $b,c$ replaced, respectively, by
\begin{align}
\label{tilde b}
\nonumber\tilde b&=b-k\left(c-1+k\right)
,\\[1ex]
 \tilde c&=c+2k.
\end{align}

Moreover the discriminant $\tilde D$ and the parameters $\tilde s_{1,2}$ of $\tilde L$ are given by
\begin{align}\label{tilde D gamma}
\tilde D&=D, \quad\tilde s_{1,2}=s_{1,2}+k.
\end{align}

Choosing $k=-s_i$, $i=1,2$, we get  $\tilde b=0$, $\tilde c_i=c-2s_i$ and therefore  $$T_{-s_i,0}^{-1}\mathcal LT_{-s_i,0}=y^{\alpha_1}\Delta_x+y^{\alpha_2} \tilde B^i_y=y^{\alpha_1}\Delta_x+y^{\alpha_2} \left(D_{yy}+\frac{c-2s_i}{y}D_y\right).$$

\begin{teo} \label{complete} Let  $\alpha_2-\alpha_1<2$ and
$$s_1+ \alpha_1^-<\frac{m+1}p<s_2+2-\alpha_2.$$   
Then $\mathcal L^{\alpha_1,\alpha_2}$ generates a bounded analytic semigroup  in $L^p_m$ which has maximal regularity. Moreover,
\begin{equation}
 \label{dominioTrasf}
D(\mathcal L^{\alpha_1,\alpha_2}_{m,p})
=y^{-s_1} W^{2,p}_{\mathcal{N}}\left(\alpha_1,\alpha_2,m-s_1p\right).
\end{equation}	

Finally, the estimate
\begin{equation} \label{splitting}
\|y^{\alpha_1}D_{x_ix_j} u\|_{L^p_m}+\|y^{\alpha_2} L_y u\|_{L^p_m} \leq C \|\mathcal L^{\alpha_1, \alpha_2} u\|_{L^p_m}
\end{equation}
holds for every $u \in D(\mathcal L^{\alpha_1,\alpha_2}_{m,p})$.
\end{teo}
{\sc Proof.}  According to the discussion above the map
$T_{-s_1,0}:L^p_{m-s_1p}\to L^p_m$
transforms  $\mathcal L^{\alpha_1,\alpha_2}$ into
$
y^{\alpha_1}\Delta_x+y^{\alpha_2}\tilde B^n_y
$ where $\tilde B^n_y=D_{yy}+\frac{\tilde c}{y}D_y ,\quad
 \tilde c=c-2s_1$. 

Since $s_1+ \alpha_1^-<\frac{m+1}p<s_2+2-\alpha_2$ is equivalent to $ \alpha_1^-<\frac{m-ps_1+1}p<\tilde c+1-\alpha_2$, 
the statement on generation and maximal regularity is therefore a translation to  $\mathcal L^{\alpha_1,\alpha_2}$ and in $L^p_m$ of the results of Section 7 for  $y^{\alpha_1}\Delta_x+y^{\alpha_2}\tilde B^n_y$ in $L^p_{m-s_1p}$. 

Also $D(\mathcal L^{\alpha_1,\alpha_2}_{m,p})=T_{-s_1,0}\left( W^{2,p}_{\mathcal{N}}\left(\alpha_1,\alpha_2,m-s_1p\right)\right)$.
Finally, \eqref{splitting} holds since the similar statement holds for $y^{\alpha_1}\Delta_x+y^{\alpha_2}\tilde B^n_y$ in $L^p_{m-s_1p}$ and 
 $$T_{-s_1,0}^{-1}\left( y^{\alpha_1} D_{x_ix_j}\right) T_{-s_1,0}=y^{\alpha_1}D_{x_ix_j}, \qquad   T_{-s_1,0}^{-1}\left(y^{\alpha_2} L_y\right)T_{-s_1,0}=y^{\alpha_2} \tilde B_y.$$
\qed

The following corollary explains why we use the term Dirichlet boundary conditions.

\begin{cor} \label{cor1}
	Let  $\alpha_2-\alpha_1<2$ and $s_1+ \alpha_1^-<\frac{m+1}p<s_2+2-\alpha_2.$ 
		\begin{itemize}
			\item[(i)] If $D>0$ then 
			\begin{equation*}
				D(\mathcal L^{\alpha_1,\alpha_2}_{m,p})=\left\{u \in  W^{2,p}_{loc}(\R^{N+1}_+): u,\ y^{\alpha_1}\Delta_xu, y^{\alpha_2}L_yu\in L^p_m\text{\;and\;}\lim_{y\to 0} y^{s_2}u(x,y)=0\right\}.
			\end{equation*}
			\item[(ii)] If $D=0$ then $s_1=s_2$ and
			\begin{equation*}
			D(\mathcal L^{\alpha_1,\alpha_2}_{m,p})=\left\{u \in  W^{2,p}_{loc}(\R^{N+1}_+): u,\ y^{\alpha_1}\Delta_xu, y^{\alpha_2}L_yu\in L^p_m\ \text{\;and\;}\lim_{y\to 0} y^{s_2}u(x,y)\in \C\right\}.
			\end{equation*}
		\end{itemize} 
		
	\end{cor}
	{\sc Proof.} Since $\tilde{c}=c-2s_1=1+2\sqrt{D}\geq 1$, both points follow by the previous theorem and \ref{trace u in W op}. \qed

\begin{os} Equality \eqref{dominioTrasf} says that $u \in D(\mathcal L^{\alpha_1,\alpha_2}_{m,p})$ if and only for every $i,j=1, \dots N$ all functions
$$u,\ y^{\frac{\alpha_1}{2}} D_{x_i}u,\  y^{\alpha_1}D_{x_i x_j}u,\ y^{\alpha_2-1}\left (D_y u+s_1\frac{u}{y} \right ),
\ y^{\alpha_2}L_y u $$
belong to $L^p_m$
but one cannot deduce, in general, that $y^{\alpha_2-1}D_y u$ and $y^{\alpha_2}D_{yy}u$ belong to $L^p_m$, as one can check on functions like $y^{-s_1}u(x)$, $u\in C_c^\infty(\R^N)$, near $y=0$. This is however possible in the special case below.
\end{os}

\begin{cor} \label{Cor-Rellich} Let  $\alpha_2-\alpha_1<2$ and
$s_1+ 2-\alpha_2<\frac{m+1}p<s_2+2-\alpha_2$.   Then $D(\mathcal L^{\alpha_1,\alpha_2}_{m,p})=W^{2,p}_{\mathcal R}(\alpha_1,\alpha_2,m)$.
\end{cor}
{\sc Proof.} Observe that  $s_1+ 2-\alpha_2>s_1+\alpha_1^- $, since  $\alpha_2<2$, $\alpha_2-\alpha_1<2$. By  Theorem \ref{complete} and Proposition \ref{RN}
 $$D(\mathcal L^{\alpha_1,\alpha_2}_{m,p})=y^{-s_1}\left( W^{2,p}_{\mathcal{N}}\left(\alpha_1,\alpha_2,m-s_1p\right)\right)=W^{2,p}_{\mathcal R}\left(\alpha_1,\alpha_2,m\right)$$
under the assumption $\frac{m-ps_1+1}p>2-\alpha_2$ which is equivalent to $s_1+ 2-\alpha_2<\frac{m+1}p$.
\qed

Concerning the mixed derivatives, we have the following result.

\begin{cor} \label{mixedcomplete} Let  $\alpha_2-\alpha_1<2$ and
$$s_1+ \alpha_1^-<\frac{m+1}p<s_2+2-\alpha_2, \quad \frac{m+1}p>s_1+1-\frac{\alpha_1+\alpha_2}{2}.$$
Then
$$
\|y^{\frac{\alpha_1+\alpha_2}{2}-1}D_{x_i} u\|_{L^p_m}+\|y^{\frac{\alpha_1+\alpha_2}{2}}D_{x_iy} u\|_{L^p_m} \leq C \|\mathcal L^{\alpha_1, \alpha_2} u\|_{L^p_m}
$$
for every $u \in D(\mathcal L^{\alpha_1,\alpha_2}_{m,p})$.
\end{cor}
{\sc Proof.} Let us write $u=y^{-s_1}v$ with $v\in W^{2,p}_{\mathcal{N}}\left(\alpha_1,\alpha_2,m-s_1p\right)$. Then 
$$
y^{\frac{\alpha_1+\alpha_2}{2}}D_{x_iy} u=y^{\frac{\alpha_1+\alpha_2}{2}}\left ( y^{-s_1}D_{x_iy} v-s_1 y^{-s_1-1}D_{x_i}v \right ).
$$
The first term on the right hand side belongs to $L^p_m (\R^{N+1}_+)$ by Theorem \ref{complete-Bessel}  and the second by Proposition  \ref{Hardy Rellich Sob}, provided $\frac{m+1}p>s_1+1-\frac{\alpha_1+\alpha_2}{2}$. This gives the estimate for $y^{\frac{\alpha_1+\alpha_2}{2}}D_{x_iy} u$. That for $y^{\frac{\alpha_1+\alpha_2}{2}-1}D_{x_i} u$ follows similarly, using Proposition \ref{Hardy Rellich Sob} again.
\qed

Observe that the condition $\frac{m+1}p>s_1+1-\frac{\alpha_1+\alpha_2}{2}$ in the previous corollary is necessary for the integrability  of the mixed derivatives of functions like $y^{-s_1}u(x)$, $u\in C_c^\infty(\R^N)$, near $y=0$.

\begin{cor} \label{uy} Let  $\alpha_2-\alpha_1<2$ and
$$s_1+ \alpha_1^-<\frac{m+1}p<s_2+2-\alpha_2, \quad \frac{m+1}p>s_1+1-\frac{\alpha_2}{2}.$$
Then
$$
\|y^{\frac{\alpha_2}{2}}D_{y} u\|_{L^p_m} \leq C\left (\|u\|_{L^p_m}+ \|\mathcal L^{\alpha_1, \alpha_2} u\|_{L^p_m}\right ).
$$
\end{cor}
{\sc Proof. }  Let us write $u=y^{-s_1}v$ with $v\in W^{2,p}_{\mathcal{N}}\left(\alpha_1,\alpha_2,m-s_1p\right)$. Then 
$$
y^{\frac{\alpha_2}{2}}D_{y} u=y^{\frac{\alpha_2}{2}}\left ( y^{-s_1}D_{y} v-s_1 y^{-s_1-1}v \right )
$$
and the thesis follows from Proposition \ref{Hardy Rellich Sob} (i).
\qed

\medskip

The above results apply also to the operator $\mathcal L=y^{\alpha_1}\Delta_x +y^{\alpha_2}B_y$, $B_y=D_{yy}+\frac{c}{y}D_y$,  when $c<1$, so that $s_1=c-1 \neq 0$, and allow to construct a realization of $\mathcal L$  different from that of Theorem \ref{complete-Bessel}. 

\begin{cor}\label{CorBdalpha}
	Let   $\alpha_2-\alpha_1<2$, $c<1$ and  and $c-1+ \alpha_1^-<\frac{m+1}p<2-\alpha_2.$ Then $\mathcal L=y^{\alpha_1}\Delta_x +y^{\alpha_2}B_y$ with domain
			\begin{equation*}
				D(\mathcal L^{\alpha_1,\alpha_2}_{m,p})=\left\{u \in  W^{2,p}_{loc}(\R^{N+1}_+): u,\ y^{\alpha_1}\Delta_xu, y^{\alpha_2}B_yu\in L^p_m\text{\;and\;}\lim_{y\to 0} u(x,y)=0\right\}.
			\end{equation*}
generates a bounded analytic semigroup  in $L^p_m$ which has maximal regularity.
\end{cor}
{\sc Proof. } This follows from Corollary \ref{cor1} (i), since $s_1=c-1$ and $s_2=0$.
\qed
Note that the generation interval  $c-1+ \alpha_1^-<\frac{m+1}p<2-\alpha_2$ under Dirichlet boundary conditions,   is larger than $ \alpha_1^-<\frac{m+1}p<c+1-\alpha_2$  given by Theorem \ref{complete-Bessel} for Neumann boundary conditions.


Let us explain what happens in Theorem \ref{complete} if we choose the second root $s_2$ instead of $s_1$. Proceeding similarly, one proves an identical result under the condition
\begin{equation} \label{nonunique}
s_2+\alpha_1^-< \frac{m+1}{p} <s_1+2-\alpha_2.
\end{equation}
However this requires the assumption $s_2<s_1+2-\alpha_2$ which is not always satisfied. When \eqref{nonunique} holds this procedure leads to a different operator, as we explain in more detail in Section 9.2.


\section{Further results, examples and applications}
\subsection{The range of contractivity}
Here we investigate when the semigroups generated  by our operators are contractive on the positive real axis. 

Let $\mathcal L=y^{\alpha_1}\Delta_x+y^{\alpha_2} B_y^n$ with   $\alpha_2<2$, $\alpha_2-\alpha_1<2$ and
$\alpha_1^{-} <\frac{m+1}p<c+1-\alpha_2$ so that the generation conditions are satisfied and  $ C_c^\infty(\R^N)\otimes \mathcal D$ is a core.

If  $I_su(x,y)=u(s^{1-\frac{(\alpha_2-\alpha_1)}{2}}x, sy)$, then 
 $I_s^{-1} \mathcal L I_s=s^{2-{\alpha_2}} \mathcal L$ and an estimate $\|e^{t \mathcal L}\| \leq e^{\omega t}$ implies  $\|e^{t \mathcal L}\| \leq 1$ (operator norms in $L^p_m$). Therefore quasi-contractivity is equivalent to contractivity.

\begin{lem} \label{contrattivita1} The operator  $\mathcal L=y^{\alpha_1}\Delta_x+y^{\alpha_2} B_y^n$ is dissipative   on $ C_c^\infty(\R^N)\otimes \mathcal D \subset L^p_m$ if and only if $y^{\alpha_2} B$ is dissipative on $\mathcal D \subset L^p_m(\R_+)$.
\end{lem}
{\sc Proof. } For $u \in  C_c^\infty(\R^N)\otimes \mathcal D$ 
\begin{align*}
-\int_{\R^{N+1}_+} (\mathcal L u) u |u|^{p-2} y^{m} \, dx\, dy&=(p-1)\int_{\R^{N+1}_+} |\nabla_x u|^2 |u|^{p-2} y^{\alpha_1+m}\, dx\, dy\\& -\int_{\R^{N+1}_+} (y^{\alpha_2} B_y u) u |u|^{p-2} y^m \, dx\, dy
\end{align*}
and the dissipativity of $\mathcal L$ follows from that of $y^{\alpha_2} B^n$. Conversely, assuming the dissipativity of $\mathcal L$, we fix $v \in \mathcal D$, $0 \neq \phi \in C_c^\infty(\R^N)$ and consider $u_n(x,y)=\phi (x/n)v(y)$. Inserting in the above identity and letting $n \to \infty$ it follows that $-\int_0^\infty (y^{\alpha_2} B_y u) u |u|^{p-2} y^m \, dy \geq 0$.
\qed

The dissipativity of $y^{\alpha_2} B^n$ will be deduced from the case $\alpha=0$, via a change of variable.

\begin{lem} \label{Hardycontractivity} The best constant in the inequality
\begin{equation} \label{c2}  \int_0^\infty u_y^2 |u|^{p-2} y^m\ dy \geq C\int_0^\infty |u|^p y^{m-2}\, dy, \quad u \in C_c^\infty (\R_+)
\end{equation}
is $C=\left (\frac{m-1}{p}\right )^2$. When $m >1$ the inequality above holds also for every $u \in \mathcal D$.
\end{lem}
{\sc Proof. } A proof that the best constant is that indicated above can be found in \cite[Proposition 8.3]{met-soba-spi-Rellich}.

When $m>1$ and $u \in \mathcal D$, let $\phi$ be a smooth cut-off functions which is equal to $0$ in $[0,1]$ and to $1$ in $[2, \infty[$. We apply the inequality above to $u_n(y)=u(y)\phi (ny)$ and get
$$
C\int_0^\infty |u_n|^p y^m\, dy \leq  \int_0^\infty u_y^2 \phi(ny)u(y)|\phi(ny)u(y)|^{p-2} y^m\, dy+ n^2 \int_{\frac 1n}^{\frac 2n} |u|^p {\phi}^2_y (ny)|\phi(ny)|^{p-2} y^m\, dy
$$
and the last term tends to 0 as $n \to \infty$, since $m>1$ and $u, \phi, \phi_y$ are bounded. One concludes by dominate convergence.
\qed
\begin{prop} \label{alpha0}
Assume $0<\frac{m+1}{p}<c+1$. The operator $B_p^n$ is dissipative in $L^p_m(\R_+)$ if and only if
\begin{itemize}
\item [(i)] $m=c$ or
\item[(ii)] $m \geq 1$ and $\frac{m-1}{p} \leq c-1$.
\end{itemize}
\end{prop}
{\sc Proof. } For $u \in \mathcal D$, $u$ constant  in $[0,a]$ we have integrating by parts
\begin{eqnarray} \label{c1}
-\int_0^\infty (Bu) u|u|^{p-2} y^m\, dy&=(p-1)\int_a^\infty u_y^2 |u|^{p-2}y^m \, dy +\frac{(1-m)(m-c)}{p}\int_a^\infty  |u|^p y^{m-2}\, dy\\
\nonumber &-\frac{m-c}{p}|u(a)|^{p}a^{m-1}
\end{eqnarray}
and (i) is immediate.

The inequality
$$
-\int_0^\infty (Bu) u|u|^{p-2} y^m\, dy=(p-1)\int_a^\infty u_y^2 |u|^{p-2}y^m \, dy +\frac{(1-m)(m-c)}{p}\int_a^\infty  |u|^p y^{m-2}\, dy\geq 0
$$
holds in $C_c^\infty (\R_+)$ if and only if $\frac{(m-1)(m-c)}{p(p-1)} \leq  \left (\frac{m-1}{p}\right )^2$ by the above lemma, which means 
\begin{equation} \label{c3}
\frac{m-1}{p}\left (c-1-\frac{m-1}{p}\right ) \geq 0. 
\end{equation}
Therefore, dissipativity can hold when $m \geq 1$ only if (ii) holds. On the other hand, if $m>1$ and (ii) holds, then letting $ a \to 0$ in \eqref{c1} we obtain 
$$
-\int_0^\infty (Bu) u|u|^{p-2} y^m\, dy=(p-1)\int_0^\infty u_y^2 |u|^{p-2}y^m \, dy +\frac{(1-m)(m-c)}{p}\int_0^\infty  |u|^p y^{m-2}\, dy
$$
which is nonnegative since
\eqref{c2} holds in $\mathcal D$, by Lemma \ref{Hardycontractivity}. Therefore (ii) is proved for $m>1$. If $m=1$ let us observe that \eqref{c1} trivially holds when $c \geq 1$. 

Finally we consider the case $m \leq 1$ and show that $B^n$ is never dissipative for $m<1$ and $c \neq m$ or for $m=1$ and $c \leq 1$, even though \eqref{c2} can hold on $C_c^\infty (\R_+)$. 

Let assume that \eqref{c3} holds, or $c-1 \leq (m-1)/p$, otherwise dissipativity fails already on $C_c^\infty (\R_+)$,  and let $u(y)=y^{-\beta}$ for $y \geq 1$ and constant in $[0,1]$. The function $u$ is not properly in $\mathcal D$ but smoothing and cutting at infinity do not make any problem. 

Assuming $(m-1)/p<\beta$ all integrals in the right hand side of \eqref{c1} converge and a straightforward computation shows that positivity is equivalent to $$\beta\left ((p-1)\beta-(m-c) \right ) \geq 0$$ for every $\beta >(m-1)/p$. However this is false for $m<1$ since the expression above is negative between 0 and $(m-c)/(p-1)$ and $(m-1)/p <0, (m-1)/p \leq (m-c)/(p-1)$.

When $m=1$ the inequality \eqref{c3} is always verified and the positivity of \eqref{c1} on $y^{-\beta}$ is equivalent to $$\beta\left ((p-1)\beta-(1-c) \right ) \geq 0$$ for $\beta >0$ which is false for small $\beta >0$ when $c<1$.
\qed
 We can now state the final contractivity result.
\begin{teo} \label{Contra}
\begin{itemize}
\item[(i)] Assume that   $\alpha_2-\alpha_1<2$ and
$$\alpha_1^{-} <\frac{m+1}p<c+1-\alpha_2.$$   Then the semigroup generated by  $y^{\alpha_1} \Delta_x+y^{\alpha_2}B^n_y$ is contractive in $L^p_m$ if and only if 
\begin{align*}
m=c-\alpha_2\qquad\text{or}\qquad 	\frac{2-\alpha_2}{p} \leq\frac{m+1}p\leq c-1+\frac{2-\alpha_2}{p}.
\end{align*}
\item[(ii)] Assume that   $\alpha_2-\alpha_1<2$ and
$$s_1+\alpha_1^{-} <\frac{m+1}p<s_2+2-\alpha_2.$$  Then the semigroup generated by  $y^{\alpha_1} \Delta_x+y^{\alpha_2}L_y$, under Dirichlet boundary conditions, is contractive in $L^p_m$ if and only if 
$$s_1+\frac{2-\alpha_2}{p}\leq \frac{m+1}{p}\leq s_2+\frac{2-\alpha_2}{p}.$$
\end{itemize}
\end{teo}
{\sc Proof.} Concerning (i), observe that by Lemma \ref{contrattivita1}  it is enough to consider $y^{\alpha_2}B^n_y$. 
According to Proposition \ref{Isometry action},  we use the isometry 
$$T_{0,-\frac{\alpha_2}{2}}:L^p_{\tilde m}\left(\R_+\right)\to L^p_m\left(\R_+\right),\quad T_{0,-\frac {\alpha_2}{2}}u(y)=\left|1-\frac{\alpha_2}{2}\right |^{\frac 1p}u(y^{1-\frac{\alpha_2}{2}}),$$
$\tilde m=\frac{m+\frac{\alpha_2}{ 2}}{1-\frac{\alpha_2}{2}}$, 
under whose action  $y^\alpha_2 B_y^n$ becomes  isometrically  equivalent to  
$
\left(1-\frac{\alpha_2}{2}\right)^2\tilde B
$ where $\tilde B=D_{yy}+\frac{\tilde c}{y}D_y $ and  
$
 \tilde c=\frac{c-\frac{\alpha_2} {2}}{1-\frac{\alpha_2} {2}}$.

 The dissipativity for $y^{\alpha_2} B$ in $L^p_m$ is then  immediate consequence of that of $\tilde{B}$ in $L^p_{\tilde m}$ already proved in Proposition \ref{alpha0}.\\

Concerning (ii), observe that, as in the previous Section, the map
$T_{-s_1,0}:L^p_{m-s_1p}\to L^p_m$
transforms  $\mathcal L^{\alpha_1,\alpha_2}$ into
$
y^{\alpha_1}\Delta_x+y^{\alpha_2}\tilde B^n_y
$ where $\tilde B^n_y=D_{yy}+\frac{\tilde c}{y}D_y ,\quad
 \tilde c=c-2s_1$. 
Therefore the dissipativity of  $\mathcal L^{\alpha_1,\alpha_2}$ in $L^p_m$ follows from that of   $y^{\alpha_1}\Delta_x+y^{\alpha_2}\tilde B^n_y$ in $L^p_{m-s_1p}$ proved in (i). 
We have that $\mathcal L^{\alpha_1,\alpha_2}$ is dissipative in $L^p_m$ if and only if $m-ps_1=c-2s_1-\alpha_2$ or $m-ps_1\geq 1-\alpha_2$ and $\frac{m-ps_1-1+\alpha_2}{p}\leq c-2s_1-1$. 
The claim follows since $m-ps_1\geq 1-\alpha_2$ and $\frac{m-ps_1-1+\alpha_2}{p}\leq c-2s_1-1$ are equivalent respectively to $\frac{m+1}{p}\geq s_1+\frac{2-\alpha_2}{p}$ and $\frac{m+1}{p}\leq s_2+\frac{2-\alpha_2}{p}$ and after observing that  $m-ps_1=c-2s_1-\alpha_2$ is equivalent to $\frac{m+1}{p}=s_1+\frac{2-\alpha_2}{p}+\frac{s_2-s_1}{p}$ and obviously $s_1+\frac{2-\alpha_2}{p}<s_1+\frac{2-\alpha_2}{p}+\frac{s_2-s_1}{p}<s_2+\frac{2-\alpha_2}{p}$.
\qed

\subsection{Further generation results and uniqueness} Let $\mathcal L=y^{\alpha_1}\Delta_x+y^{\alpha_2} L_y$, $\alpha_2<2$, and keep the notation of Section 8, in particular $\mathcal L^{\alpha_1, \alpha_2}_{m,p}$ is the operator constructed therein. Let us define the maximal operator $\mathcal L^{max}_{m,p}$ as $\mathcal L$ on the maximal domain
\begin{equation} \label{Dmax}
D(\mathcal L^{max}_{m,p})=\{u \in L^p_m\cap W^{2,p}_{loc}(\R^{N+1}_+): \mathcal Lu \in L^p_m\}
\end{equation}
and the minimal operator  $\mathcal L^{min}_{m,p}$ as the closure of $\mathcal L$ initially defined on $C_c^\infty (\R^{N+1}_+)$. By local elliptic regularity  $\mathcal L^{max}_{m,p}$ is closed and then, since  $(\mathcal L, C_c^\infty (\R^{N+1}_+))$ admits the closed extension $\mathcal L^{max}_{m,p}$, its closure is well defined. Clearly
$\mathcal L^{min}_{m,p} \subset \mathcal L^{\alpha_1, \alpha_2}_{m,p} \subset \mathcal L^{max}_{m,p}$.

Integrating by parts one sees that the formal adjoint of $\mathcal L$ is the operator $\mathcal L^*=y^{\alpha_1}\Delta_x+y^{\alpha_2} L^*_y$ in $L^{p'}_m$ where
$$L^*_y=D_{yy}+\frac{\tilde c}{y}D_y-\frac{\tilde b}{y^2}, \quad \tilde c=2\alpha_2+2m-c, \quad \tilde b=b-(\alpha_2+m-c)(\alpha_2+m-1).$$ Moreover, the characteristic numbers of $\mathcal L^*$ are given by $$D^*=D,\quad  s_1^*=\alpha_2+m-1-s_2,\quad  s_2^*=\alpha_2+m-1-s_1.$$

\begin{lem} \label{dual} The dual of  $\mathcal L^{min}_{m,p}$ is  $\mathcal L^{*,max}_{m,p'}$ and the dual of  $\mathcal L^{max}_{m,p}$ is   $\mathcal L^{*,min}_{m,p'}$.
\end{lem}
{\sc Proof.} Since $L^p_m$ is reflexive, it is sufficient to prove the first equality. The second follows by duality from the first, changing $p$ with $p'$.
If $u \in C_c^\infty (\R^{N+1}_+)$ and $v \in D(\mathcal L^{*,max}_{m,p'})$ one can integrate by parts and get
$$
\int_{\R^{N+1}_+}v(\mathcal L u)\, y^m \, dx\, dy=\int_{\R^{N+1}_+}u( \mathcal L^* v )\,  y^m \, dx\, dy
$$
and hence  $\mathcal L^{*,max}_{m,p'}$ is a restriction of the dual of $\mathcal L^{min}_{m,p}$. Conversely, if $v \in L^{p'}_m$ and
 $$
\int_{\R^{N+1}_+}v(\mathcal L u)\, y^m \, dx\, dy=\int_{\R^{N+1}_+}u f   y^m \, dx\, dy
$$
for some $f \in L^{p'}_m$, then by local elliptic regularity (the coefficients of $\mathcal L$ are smooth in the interior of $\R^{N+1}_+$), $v \in W^{2,p'}_{loc}(\R^{N+1}_+)$ and ${\mathcal L}^*v=f \in L^{p'}_m$.
\qed

\begin{prop} \label{Lmin} Let  $\alpha_2-\alpha_1<2$ and
$s_1+ 2-\alpha_2 \leq \frac{m+1}p<s_2+2-\alpha_2$.   Then $\mathcal L^{\alpha_1, \alpha_2}_{m,p}=\mathcal L^{min}_{m,p}$.
\end{prop}
{\sc Proof.} Observe that  $s_1+ 2-\alpha_2>s_1+\alpha_1^- $, since  $\alpha_2<2$, $\alpha_2-\alpha_1<2$. By  Theorem \ref{complete} 
 $$D(\mathcal L^{\alpha_1,\alpha_2}_{m,p})=y^{-s_1}\left( W^{2,p}_{\mathcal{N}}\left(\alpha_1,\alpha_2,m-s_1p\right)\right).$$
The assumption $s_1+ 2-\alpha_2 \leq \frac{m+1}p$ is equivalent to $\frac{m-ps_1+1}p\geq 2-\alpha_2$ and one concludes by \ref{Core C c infty}.
\qed

Note that when $s_1+ 2-\alpha_2 < \frac{m+1}p<s_2+2-\alpha_2$, then we have also $D(\mathcal L^{\alpha_1,\alpha_2}_{m,p})=W^{2,p}_{\mathcal{R}}(\alpha_1,\alpha_2,m)$, by Corollary \ref{Cor-Rellich}.

\begin{prop} \label{Lmax} Let  $\alpha_2-\alpha_1<2$ and
$s_1< \frac{m+1}p \leq s_2$.   Then $\mathcal L^{max}_{m,p}$ generates an analytic semigroup.
\end{prop}
{\sc Proof.} Let us  consider the adjoint $\mathcal L^*$. Then $s_1^*+ 2-\alpha_2 \leq \frac{m+1}{p'}<s_2^*+2-\alpha_2$ and then, by the proposition above,  $\mathcal L^{*,min}_{m,p'}$ generates a semigroup in $L^{p'}_m$. By standard semigroup duality in reflexive spaces, the  dual operator  $\mathcal L^{max}_{m,p}$, see Lemma \ref{dual}, generates a semigroup in $L^p_m$.
\qed

Observe that when $\alpha_1<0$, $\mathcal L^{max}_{m,p}$ generates a semigroup also when the condition $s_1+\alpha_1^- <\frac{m+1}{p}$ is violated. However, if this last holds, then $\mathcal L^{\alpha_1, \alpha_2}_{m,p}=\mathcal L^{max}_{m,p}$.

\begin{prop} Let  $\alpha_2-\alpha_1<2$ and
$s_1+ \alpha_1^-< \frac{m+1}p \leq s_2$.   Then $\mathcal L^{\alpha_1, \alpha_2}_{m,p}=\mathcal L^{max}_{m,p}$.
\end{prop}
{\sc Proof.} In fact  $\mathcal L^{\alpha_1, \alpha_2}_{m,p}$ is well defined and generates a semigroup. Since  $\mathcal L^{max}_{m,p}$ extends $\mathcal L^{\alpha_1, \alpha_2}_{m,p}$ and both are generators, they coincide. 
\qed

By duality, we can extend the generation interval.

\begin{prop} \label{rangegeneration} If  $\alpha_2-\alpha_1<2$ and $\frac{m+1}{p} \in (s_1, s_2+2-\alpha_2-\alpha_1^-) \cup (s_1+\alpha_1^-, s_2+2-\alpha_2)$
 a realization  $\mathcal L^{min}_{m,p} \subset \mathcal L_D \subset \mathcal L^{max}_{m,p}$ generates a semigroup.
\end{prop}
{\sc Proof.} In fact, if $s_1+\alpha_1^-< \frac{m+1}p  <s_2+2-\alpha_2$ we can take  $\mathcal L^{\alpha_1, \alpha_2}_{m,p}$ and if $s_1< \frac{m+1}p < s_2+2-\alpha_2-\alpha_1^-$ we can take the adjoint of $\mathcal{ L}^{*, \alpha_1, \alpha_2}_{m,p'}$ since the condition is equivalent to $s^*_1+\alpha_1^-< \frac{m+1}{p'}  <s_2^*+2-\alpha_2$.
\qed

For the $1d$ operator $y^{\alpha_2}L_y$ it is known that a such a realization exists if and only if $s_1<\frac{m+1}{p} <s_2+2-\alpha_2$, see \cite[Theorem 1.1, Theorem 1.2]{MOSS} for the case $m=0$.  For general $m$ and  $\alpha_2=0$ see \cite[Propositions 2.4, 2.5]{MNS-Max-Reg} from which, by the transformation $T_{0,-\frac{\alpha_2}{2}}$, it is possible deduce the general case. However the above proposition yields a semigroup in this interval only when $s_1+\alpha_1^- <s_2+2-\alpha_2-\alpha_1^-$. For example, when $\alpha_1=\alpha_2=\alpha<0$ this requires $|\alpha|<s_2-s_1+2$.

Let us show that when the condition $s_1+\alpha_1^- < \frac{m+1}{p}$ is violated, the regularity estimate 
\begin{equation} \label{reg}
\|y^{\alpha_1}D_{x_ix_j} u\|_{L^p_m}+\|y^{\alpha_2} L_y u\|_{L^p_m} \leq C \|\mathcal L u\|_{L^p_m}
\end{equation}
may fail for $u$ in the domain of the operator.

\begin{esem} \label{noreg}  Let $\mathcal L=y^{-\beta}(\Delta_x+D_{yy})$, $\beta>0$. Then $s_1=-1$, $s_2=0$, $\mathcal L$ generates under Neumann  boundary conditions when $\beta<\frac{m+1}{p}<1+\beta$ and under Dirichlet boundary conditions when $-1+\beta<\frac{m+1}{p}<2+\beta$ and both operators satisfy \eqref{reg}. 

However, when $-1<\frac{m+1}{p} \leq (-1+\beta)\wedge 0$, $\mathcal L^{max}_{m,p}$ is a generator for which \eqref{reg} fails. Indeed, let $\eta$ be a smooth function equal to $1$ in $[0,\frac 12]$ and to $0$ in $[1, \infty[$ and 
$$u(x,y)=\eta(y) \left (\frac{y+1}{\left(|x|^2+(y+1)^2\right)^{\frac{N+1}{2}}} +\frac{y-1}{\left(|x|^2+(y-1)^2\right)^{\frac{N+1}{2}}}\right ).
$$ Note that $u$ is, for small $y$, the difference of the Poisson kernels on the hyperplanes $y=\pm 1$. Then $u \in L^p_m$, since $(m+1)/p+1>0$,  and $\mathcal L u \in L^p_m$, since $\Delta u=0$ for $0 \leq y \leq \frac 12$, so that $u \in D(\mathcal L^{max}_{m,p})$. However, $y^{-\beta}\Delta_x u$ and $y^{-\beta}D_{yy}u$ do not belong to $L^p_m$, since $(m+1)/p \leq -1+\beta$. 
\end{esem}

\medskip

A natural question arises   if different boundary conditions can be imposed to produce different semigroups in $L^p_m$. This is the case, for example, for the  operator 
$\mathcal L=y^{\alpha_1}\Delta_x +y^{\alpha_2}B_y$ in Theorem \ref{complete-Bessel} and  Corollary \ref{CorBdalpha}, in the range $c<1$, $\alpha_1^-<\frac{m+1}p<c+1-\alpha_2$, where both    boundary conditions $\lim_{y\to 0}u=0$ and  $\lim_{y\to 0}y^cD_yu=0$  can be imposed and produce different semigroups.  

As in \cite[Section 5]{MNS-Caffarelli}  we look for realizations $\mathcal L_D$ such that $\mathcal L^{min}_{m,p} \subset \mathcal L_D \subset \mathcal L^{min}_{m,p}$ . From Propositions \ref{Lmin} and \ref{Lmax} it follows that $\mathcal L_D$ is unique in $L^p_m$  if $s_1 < \frac{m+1}{p} \leq s_2$ or $s_1+2-\alpha-2 \leq \frac{m+1}{p} <s_2+2-\alpha_2$. Uniqueness then  holds in the generation range of  $L_y$, namely $(s_1, s_2+2-\alpha_2)$, when these two intervals overlap, that is when 
$s_1+2-\alpha_2 \leq s_2$ or equivalently $D\geq \left(1-\frac{\alpha_2}{2}\right)^2$. In this case, uniqueness does not depend on $p$ and $m$.

Uniqueness may fail if
$ s_2<s_1+2-\alpha_2 $ and $ (m+1)/p \in (s_2,s_1+2-\alpha_2)$, as we show under the stronger assumptions $s_1 \neq s_2$ and $s_2+\alpha_1^-<\frac{m+1}p <s_1+2-\alpha_2$.
\begin{prop}
	If  $0<D<\left(1-\frac{\alpha_2}{2}\right)^2$ and $s_2+\alpha_1^-<\frac{m+1}p <s_1+2-\alpha_2$,    the  operator $
	\mathcal{L}=y^{\alpha_1}\Delta_x+y^{\alpha_2}L_y
	$ with domain
	\begin{align} \label{L neumann}
		D(\mathcal L)
		&=y^{-s_2} W^{2,p}_{\mathcal{N}}\left(\alpha_1,\alpha_2,m-s_2p\right)\\[1ex]\nonumber
		&=\left\{u \in  W^{2,p}_{loc}(\R^{N+1}_+): u,\ y^{\alpha_1}\Delta_xu, y^{\alpha_2}L_yu\in L^p_m\text{\;and\;}\lim_{y\to 0} y^{s_1+1}\left(D_yu+s_2 \frac{u}y\right)=0\right\}.
	\end{align}
generates a semigroup  in $L^p_m$.
\end{prop}
{\sc Proof.} We proceed as  in the proof of  Theorem \ref{complete} but in place of the isometry $T_{-s_1,0}$  we use  $T_{-s_2,0}:L^p_{m-s_2p}\to L^p_m$
which transform $\mathcal{L}$ into 
 $y^{\alpha_1}\Delta_x+ y^{\alpha_2}\tilde B_y^n$ where $\tilde B_y^n=D_{yy}+\frac{\tilde c}yD_y$, $\tilde c=c-2s_2$. Observe that, under the given hypotheses,   $\tilde c=1-2\sqrt D>-1+\alpha_2$ and the claim follows by Theorem \ref{complete-Bessel}.\qed

We point out that in the range $s_2+\alpha_1^-<\frac{m+1}p <s_1+2-\alpha_2$ the operators $\mathcal L^{\alpha_1,\alpha_2}_{m,p}$ of Theorem \ref{complete} and $(\mathcal L,D(\mathcal L))$ just constructed are different. In fact let $f=a(x)b(y)\in  C_c^\infty (\R^{N})\times\mathcal D$ a function in the core defined in \eqref{defD}. Then $u=y^{-s_1}f$ belongs to $D(\mathcal L^{\alpha_1,\alpha_2}_{m,p})$ but not to $D(\mathcal L)$
since $\lim_{y\to 0}y^{s_1+1}\left(D_yu+s_2 \frac{u}y\right)=s_2-s_1>0$.


		
\subsection{Baouendi-Grushin operator}
Our results apply to  generalized Baouendi-Grushin operators $$\mathcal L=y^\alpha\Delta_x+L_y, \quad \alpha>-2$$ in the half space $\R^{N+1}_+$ both with Neumann and Dirichlet boundary conditions, but 
we restrict ourselves to the classical case $\mathcal L=y^\alpha\Delta_x+D_{yy}$ in the whole space $\R^{N+1}$ with the Lebesgue measure. Our results  improve those from \cite{MNS-Grushin}, allowing negative $\alpha$ and showing maximal regularity, besides domain characterization.

\begin{prop}
Let $\alpha>-\frac{1}{p}$. Then  $\mathcal L=|y|^\alpha\Delta_x+D_{yy}$ with domain
$$D(\mathcal L)=\left\{u\in W^{2,p}_{loc}(\R^{N+1}):\ u,\  y^{\alpha} D_{x_ix_j}u,\ y^\frac{\alpha}{2} D_{x_i}u,\  D_{yy}u,\ D_{y}u,\ y^\frac{\alpha}{2} D_{x_iy}u\in L^p(\R^{N+1})\right\}$$
 generates a bounded analytic semigroup  in $L^p(\R^{N+1})$ which has maximal regularity. 
\end{prop}
{\sc Proof.} By Theorems \ref{complete-Bessel}, \ref{complete}
the  operator $\mathcal L$ generates an analytic semigroup in $L^p(\R^{N+1}_+)$ both with Dirichlet and Neumann boundary conditions. We can therefore consider the operators $\mathcal L_i$, $i=1,2$, that is $\mathcal L$ with domains 
\begin{align*}
D_1=&\Bigl\{u \in  W^{2,p}_{loc}(\R^{N+1}_+): u,\ y^{\alpha} D_{x_ix_j}u,\ y^\frac{\alpha}{2} D_{x_i}u,  D_{yy}u,\ D_{y}u, y^\frac{\alpha}{2} D_{x_iy}u\in L^p(\R^{N+1}_+), \\& \hspace{70ex}\lim_{y\to 0}u(x,y)=0\Bigr\},
\end{align*}
\begin{align*}
D_2&=W^{2,p}_{\mathcal{N}}\left(\alpha,0,0\right)\\[1ex]
&=\Bigl\{u\in W^{2,p}_{loc}(\R^{N+1}):\ u,\  y^{\alpha} D_{x_ix_j}u,\ y^\frac{\alpha}{2} D_{x_i}u,  D_{yy}u,\ D_{y}u,y^\frac{\alpha}{2} D_{x_iy}u\in L^p(\R^{N+1}_+), \\& \hspace{70ex}\lim_{y\to 0}D_yu(x,y)=0\Bigr \}.
\end{align*}
The mixed derivatives estimates follows from Theorem \ref{complete-Bessel} and Corollary \ref{mixedcomplete}.

Let  $P_1, P_2: L^p(\R^{N+1})\to L^p(\R^{N+1}_+)$ be the even and odd projections
$$(P_1f)(x,y)=\frac{f(x,y)+f(x,-y)}{2},\quad (P_2f)(x,y)=\frac{f(x,y)-f(x,-y)}{2}$$ and 
$E_1, E_2: L^p(\R^{N+1}_+) \to : L^p(\R^{N+1})$ the even and odd extensions
\begin{align*}
E_1u(x,y)= \begin{cases}
u(x,y),\qquad &\text{if}\quad y>0;\\[1ex]
u(x,-y),\qquad &\text{if}\quad y<0;\\[1ex]
\end{cases}
\end{align*}
\begin{align*}
E_2u(x,y)= \begin{cases}
u(x,y),\qquad &\text{if}\quad y>0;\\[1ex]
-u(x,-y),\qquad &\text{if}\quad y<0.\\[1ex]
\end{cases}
\end{align*}
Note that $E_1P_1+E_2P_2=I_{L^p(\R^{N+1})}$, $P_i (D(\mathcal L)) \subset D_i$, $E_i (D_i) \subset D(\mathcal L)$, $i=1,2$ and that $D(\mathcal L)=D_1 \oplus D_2$  algebraically and topologically with respect to the Sobolev norm. Then $\mathcal L=E_1 \mathcal L_1P_1+E_2\mathcal L_2 P_2$ and everything follows from the properties of $\mathcal L_1, \mathcal L_2$.
\qed

\subsection{The operator $y\Delta_x+yB_y$} \label{alpha=1} 

We specialize  and comment here the results obtained in the special case $\alpha=1$, that is  for $\mathcal L=y\Delta_x+yB_y=y\Delta_x+yD_{yy}+cD_y$ where $B_y=D_{yy}+\frac{c}{y}D_y$.

Theorem \ref{complete-Bessel} applies when $0<\frac{m+1}{p} <c$ and yields generation and all other properties listed therein for $y\Delta_x+yB^n_y$ with domain 
$$W^{2,p}_{\mathcal N}(1,1,m)=\{u \in W^{2,p}_{loc}(\R^{N+1}_+): u,\, \nabla_x u,\, D_y u,\, yD_{x_ix_j}u,\, yD_{yy}u,\, yD_{x_iy}u  \in L^p_m\},$$
using also Proposition \ref{Hardy Rellich Sob} for $\nabla_x u$.

This result has been already proved in \cite{Koch}  but also in \cite{Pruss-Deg} when $c \geq 1$ and $m=0$ and in \cite{FornMetPallScn3} when $p=2, m=0$ (and $c>\frac 12$).

Note that, when $m=0$, then $$W^{2,p}_{\mathcal N}(1,1,0)=\{u \in W^{1,p}(\R^{N+1}_+):  yD_{x_ix_j}u,\, yD_{yy}u,\, yD_{x_iy}u  \in L^p(\R^{N+1}_+)\}$$
and the associated elliptic and parabolic problems seem to have no boundary condition. In our approach, the Neumann boundary condition is indeed imposed to $yD_y u$, by requiring that 
$\frac 1y (yD_yu) \in L^p(\R^{N+1}_+)$.

Theorem \ref{complete} says nothing new when $c \geq 1$, since then $s_1=0$ and the transformation $T_{-s_1,0}$ is the identity. However, when $c<1$ then  $s_1=c-1, s_2=0$ and Theorem \ref{complete} yields a different operator $\mathcal L^{1,1}_{m,p}$ in the range $c-1<\frac{m+1}{p}<1$. Its domain is
\begin{equation*}
				\left\{u \in  W^{2,p}_{loc}(\R^{N+1}_+): u,\ yD_{x_ix_j}u,\, yD_{x_iy}u,\,  yD_{yy}u+cD_yu \in L^p_m\text{\;and\;}\lim_{y\to 0}u(x,y)=0\right\}
			\end{equation*}
by Corollary \ref{cor1}(i) and Corollary \ref{mixedcomplete} for the mixed derivative. However, it is not true that $yD_{yy}u$ and $D_y u$ belong to $L^p_m$ separately, even when $m=0$, see also \cite{FornMetPallScn4}. On the other hand,  when $c<\frac{m+1}{p}<1$, then Corollary \ref{Cor-Rellich} applies and gives 
$$
D(\mathcal L^{1,1}_{m,p})=W^{2,p}_{\mathcal R}(1,1,m).
$$
In particular, if $m=0$, it follows that for $c<\frac 1p$ $$D(\mathcal L^{1,1}_{0,p})=\{u \in W^{1,p}_0 (\R^{N+1}_+): yD_{x_ix_j}u,\, yD_{yy}u,\, yD_{x_iy}u  \in L^p(\R^{N+1}_+)\},$$  a result already proved in \cite{FornMetPallPruss}.

Finally, let us specialize the results of  Section 9.2, see also \cite{FornMetPallScn2}. If $c<1$
\begin{itemize}
	\item[(i)] $\mathcal L^{1,1}_{m,p}=\mathcal L^{min}_{m,p}$ when $c \leq (m+1)/p < 1$;\quad  $\mathcal L^{1, 1}_{m,p}=\mathcal L^{max}_{m,p}$ when $c-1<(m+1)/p\leq 0$;
	\item[(ii)] uniqueness fails if and only if $0<(m+1)/p< c$.
\end{itemize}
 Instead, if $c \geq 1$
 \begin{itemize}
 	\item[(i)] $\mathcal L^{1, 1}_{m,p}=\mathcal L^{min}_{m,p}$ when $1 \leq (m+1)/p < c$;\quad $\mathcal L^{1, 1}_{m,p}=\mathcal L^{max}_{m,p}$ when $0<(m+1)/p\leq c-1$;
 	\item[(ii)]uniqueness fails if and only if $c-1<(m+1)/p < 1$.
 \end{itemize} 

\bibliography{../TexBibliografiaUnica/References}

\end{document}